\documentclass[12pt,a4paper]{article}
\usepackage[cm,myheadings]{fullpage}
\usepackage{amssymb,amsmath}
\usepackage{rotating,booktabs}
\usepackage{colordvi}
\makeatletter
\def\@oddhead{\ifnum\thepage>1\vbox{\hbox to \textwidth{%
                 \ifodd\thepage\relax%
		 \textbf{\thepage}%
		 \hfil\textsc{Gordon James and Andrew Mathas}\hfil%
		 \else%
		   \hfil\textsc{Symmetric group blocks of small defect}%
		   \hfil\textbf{\thepage}%
		  \fi}}%
	      \fi}
\let\bigcirc=\bullet
\let\bar=\overline
\let\frak\mathfrak

\newbox\squ  % box character for ends of proofs
\setbox\squ=\hbox{\vrule width.6pt
            \vbox{\hrule height.6pt width.4em\kern1ex\hrule height.6pt}%
            \vrule width.6pt}
\def\endproof{%
    \ifmmode\eqno\copy\squ\medskip\else{\unskip\nobreak\hfil%
    \penalty50\hskip2em\hbox{}\nobreak\hfil\copy\squ
    \parfillskip=0pt \finalhyphendemerits=0\penalty-100\medskip}\fi
}
\begin{document}
%\pagewiselinenumbers

\[
\text{\textbf{SYMMETRIC\ GROUP\ BLOCKS\ OF\ SMALL\ DEFECT}} 
\]

\[
\text{Gordon James}^{a}\text{ and Andrew Mathas}^{b} 
\]
\[
^{a}\text{\textit{Department of Mathematics, Imperial College, 180 Queen's
Gate, London SW7 2BZ, UK}} 
\]

\[
^{b}\text{\textit{School of Mathematics F07, University of Sydney, Sydney
NSW 2006, Australia}} 
\]

\bigskip 
\[
\text{1. INTRODUCTION} 
\]

This paper is concerned with the modular representation theory of
the symmetric groups. Throughout, we fix a positive integer $n$ and a
prime $p>0$ and we consider representations of the symmetric group
$\frak{S}_{n}$ of degree $n$ over a field of characteristic $p$. We
adopt the standard notation for the representation theory of the
symmetric groups from [\textbf{JK}].

It is well known that a $p$-block of a symmetric group $\frak{S}_{n}$ is
determined by its $p$-core and its weight, and that the weight of a block is
equal to the defect of the block if $p$ exceeds the weight [\textbf{JK}]. In
this paper we shall be concerned mainly with blocks of small defect.

Let $\lambda$ and $\mu$ be partitions of $n$ with $\mu$ being 
$p$-regular. As usual, the symmetric group $\frak{S}_{n}$ has a Specht module 
$S(\lambda)$ and a $p$-modular irreducible module $D(\mu)$. The
decomposition number $[S(\lambda):D(\mu)]$ is defined to be the
composition multiplicity of $D(\mu)$ in $S(\lambda)$. \ The following
facts are known about blocks of weight $w$.

(a) If $w=0$ or $1$ then all the decomposition numbers of the block are $0$
or $1$.

(b) If $w=2$ and $p>2$ then all the decomposition numbers of the block are 
$0$ or $1$ [\textbf{S2}].

(c) If $w=4$ then some decomposition numbers of a block can be greater than 
$1$, even if $p>w$.

Moreover, if $w=0,1$ or $2$ then there is a known method for determining the
decomposition numbers [\textbf{R}]. The situation for the case $w=3$ is
still not properly understood. In particular, the decomposition number

\[
[S(2p-2,2p-2,p-1,1):D(3p-3,2p-1)] 
\]
is yet to be determined for $p>5$. This is just one of a collection of
decomposition numbers for weight $3$ which we are unable to evaluate.

Our investigation of blocks of weight $3$ grew out of an attempt to improve
upon the earlier results of Martin and Russell [\textbf{MR}] by explicitly
calculating the decomposition numbers. This led us to discover various
errors and omissions in [\textbf{MR}] which place in doubt the claim made
there that when $p>w$ all the decomposition numbers are $0$ or $1$. Note,
incidentally, that if the decomposition numbers for a given block are known
to be $0$ or $1$ then the decomposition numbers can, in principle, be
determined by applying Schaper's Theorem [\textbf{Sch}, \textbf{JM1}].

\bigskip

\[
\text{2. BASIC\ RESULTS} 
\]

By using a result of Brundan and Kleshchev, we are able to improve upon the
presentation of several of the basic techniques used in [\textbf{S1}, \textbf{S2}, 
\textbf{MR}] for estimating decomposition numbers. 
In order to state these results recall that the \textbf{diagram} of
a partition $\lambda$ is the set of \textbf{nodes}
\[
[\lambda] = \{\, (i,j)\,\mid\,1\le j\le\lambda_i\,\}.
\]
We think of $[\lambda]$ as being an array of crosses in the plane and
we will refer to the rows and columns of $[\lambda]$ which should be
interpreted in the obvious way.

A node $x\in[\lambda]$ is \textbf{removable} if
$[\lambda]\setminus\{x\}$ is the diagram of a partition. Similarly, a
node $y\notin[\lambda]$ is \textbf{addable} if $[\lambda]\cup\{y\}$ is
the diagram of a partition. The node
$x=(i,j)$ is called an $r$-node if $r\equiv j-i\pmod p$. A
removable $r$-node $x\in[\lambda]$ is \textbf{normal} if whenever $y$ is an addable
$r$-node in $[\lambda]$ which is in an earlier row than $x$ then there are more
removable $r$-nodes between  $x$ and $y$ than there are
addable $r$-nodes [\textbf{K}].

Finally, recall that a partition $\mu$ is \textbf{$p$-regular} if
no $p$ non-zero parts of $\mu$ are equal. Then $D(\mu)\ne0$ if and
only if $\mu$ is $p$-regular.

\bigskip
\noindent 2.1 PROPOSITION \textit{Assume that }$\lambda$\textit{\ and }$\mu$\textit{\
are partitions of }$n$\textit{\ with }$\mu$\textit{\ being $p$-regular, and that 
}$k$ \textit{is a positive integer such that}
\begin{enumerate}
    \item $\lambda$\textit{\ has at most }$k$\textit{\ removable $r$-nodes; and}
    \item $\mu$\textit{\ has a least }$k$\textit{\ normal $r$-nodes.}
\end{enumerate}
\noindent %
\textit{Then $[S(\lambda):D(\mu)]$ is either zero or is
equal to an explicit decomposition number of
$\frak{S}_{n-k}$. More precisely, if $\lambda$ has fewer than $k$
removable $r$-nodes then $[S(\lambda):D(\mu)]=0$; and if 
$\lambda$ has exactly $k$ removable $r$-nodes
then $[S(\lambda):D(\mu)]=[S(\bar{\lambda}):D(\bar{\mu})]$,
where $\bar{\lambda}$ is the partition obtained from $\lambda$ by removing its $k$
$r$-nodes, and $\bar{\mu}$ is the partition obtained from $\mu$ by
removing its lowest $k$ normal $r$-nodes.
}

\bigskip

\noindent\textit{Proof.} By $r$-restricting $D(\mu)$ $k$ times we obtain an 
$\frak{S}_{n-k}$-module which contains $D(\bar{\mu})$ as a submodule. If $\lambda$
has fewer than $k$ removable $r$-nodes then by $r$-restricting 
$S(\lambda)$ $k$ times we obtain the zero module, so $[S(\lambda):D(\mu)]=0$. If 
$\lambda$ has exactly $k$ removable $r$-nodes then 
$[S(\lambda):D(\mu)]=[S(\bar{\lambda}):D(\bar{\mu})]$ by [\textbf{BK}, Lemma 2.13].
\endproof

We now recall the notion of an \textbf{abacus} from [\textbf{J2}]. A
$p$-abacus has $p$ runners, which we label as runner $1$ to runner
$p$, reading from left to right. The bead positions on the abacus are
labelled $1,2,3,\dots$, reading from left to right and then top to
bottom.  Thus, the beads on runner $r$ have  labels $r+pk$, for some
$k\ge0$.  

Recall that if $\lambda=(\lambda_1,\lambda_2,\dots)$ with
$\lambda_i=0$, whenever $i>k$, then $\lambda$ has an abacus
configuration with $k$ beads at positions $\{\lambda_i+k-i+1\,\mid1\le
i\le k\}$.  Note that if $\lambda$ has $k$ non-zero parts then
$\lambda$ can be represented on an abacus with $k'$ beads whenever
$k'\ge k$. For example the partition $(15,13,6,4^2,2^2)$ can be
represented as an abacus with $10$ and $11$ beads, respectively, as
follows:
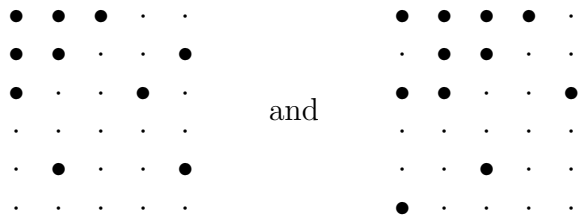
\[
\begin{array}{ccccc}
\bigcirc & \bigcirc & \bigcirc & \cdot & \cdot \\ 
\bigcirc & \bigcirc & \cdot & \cdot & \bigcirc \\ 
\bigcirc & \cdot & \cdot & \bigcirc & \cdot \\ 
\cdot & \cdot & \cdot & \cdot & \cdot \\ 
\cdot & \bigcirc & \cdot & \cdot & \bigcirc \\ 
\cdot & \cdot & \cdot & \cdot & \cdot
\end{array}
\text{\qquad and\qquad}
\begin{array}{ccccc}
\bigcirc & \bigcirc & \bigcirc & \bigcirc & \cdot \\ 
\cdot & \bigcirc & \bigcirc & \cdot & \cdot \\ 
\bigcirc & \bigcirc & \cdot & \cdot & \bigcirc \\ 
\cdot & \cdot & \cdot & \cdot & \cdot \\ 
\cdot & \cdot & \bigcirc & \cdot & \cdot \\ 
\bigcirc & \cdot & \cdot & \cdot & \cdot
\end{array}
\]
 An abacus representation with $k$ beads can be converted into one with
$k+1$ beads by shifting all beads one position to the right and then
adding a new bead at position $1$.
 
It is convenient to say that a bead on runner $r$ is an $r$-node. This
changes the definition of $r$-node above by a constant and causes
no harm. With this convention, removing an $r$-node from a partition
corresponds to moving a bead on runner $r$ one space to the left (with
an obvious modification if $r=1$), and adding an $r$-node corresponds
to moving a bead on runner $r-1$ one space to the right (with an
obvious modification if $r=1$).

By definition, a \textbf{$p$-core} is partition which has an
abacus configuration in which all of the beads are positioned as high
as possible  on each runner. A partition has \textbf{$p$-weight}
$w$ if its abacus configuration can be obtained by starting with
the abacus configuration of a $p$-core and sliding $w$ (not
necessarily distinct) beads down one  position on their runner. In this
way, we attach a $p$-core to each partition of weight $w$.

Finally, recall that all of the irreducible constituents of a Specht
module $S(\lambda)$ belong to the same block and, further, that
$S(\lambda)$ and $S(\mu)$ belong to the same block if and only if
$\lambda$ and $\mu$ have the same $p$-core [\textbf{JK}]. Consequently,
$S(\lambda)$ and $S(\mu)$ belong to the same block if and only if they
are of the same weight and they have abacus configurations which have
the same number of beads on each runner. We will say that two
partitions $\lambda$ and $\mu$ belong to a block $B$ if $S(\lambda)$ and
$S(\mu)$ are both contained in~$B$.

We can now present some corollaries of Proposition 2.1.

\bigskip

\noindent 2.2 COROLLARY \textit{Suppose that the partition $\lambda$ of $n$
has exactly $k$ removable $r$-nodes and no
addable $r$-nodes.  Let $\mu$ be a $p$-regular partition of $n$. Then 
$[S(\lambda):D(\mu)]$ is equal to an explicit
decomposition number of $\frak{S}_{n-k}$ which is 
in a block of the same weight as $\lambda$.
}

\bigskip

\noindent\textit{Proof.}\textbf{\ }We may assume that $\mu$ is in the same block as 
$\lambda$. Hence $\mu$ has exactly $k$ more removable $r$-nodes than
addable $r$-nodes and so has at least $k$ normal $r$-nodes. The Corollary
now follows immediately from Proposition 2.1. (The remark that the block of 
$\frak{S}_{n-k}$ has the same weight as $\lambda$ follows from the fact that
the abacus configuration of $\bar{\lambda}$ can be obtained from that of
$\lambda$ by swapping runners $r-1$ and $r$.)
\endproof

\noindent 2.3 COROLLARY \textit{Suppose that $B$ is a block
of $\frak S_n$ with the property that for every partition in $B$ there
exists an $r$ such that the partition has a removable $r$-node but no
addable $r$-node.  Then we can equate each decomposition number of~$B$
with an explicit decomposition number for a smaller symmetric group.}

\bigskip

From now on, we assume that we are dealing with a block of weight $w$.

\bigskip

\noindent 2.4 COROLLARY \textit{Suppose that $w\leq 3$. Then every
decomposition number for the principal block of $\frak{S}_{wp}$ is
either zero or can be equated with an explicit decomposition number of
$\frak{S}_{wp-1}$.}

\bigskip

\noindent\textit{Proof.} The $p$-core of the principal block of
$\frak{S}_{wp}$ is empty, and so it can be represented on an abacus
with $w$ beads on each runner, with all the beads pushed as far up as
possible. Suppose that $S(\lambda)$ belongs to the principal block of
$\frak{S}_{wp}$, so that the abacus configuration for $\lambda$ is
obtained from the $p$-core configuration by moving  $w$ beads, not
necessarily distinct, down one position on their runners.  Since
$w\leq 3$, we see that in the abacus configuration for $\lambda$, for
each $r$, we can move at most one bead from runner $r$ to runner
$r-1$. In other words, $\lambda$ has at most one removable $r$-node.
Now suppose that $\mu$ is $p$-regular and choose a normal  $r$-node
of $\mu$, for some $r$.  Proposition 2.1 now allows us to deduce that
$[S(\lambda):D(\mu)]$ is either zero or equal to a decomposition
number of $\frak{S}_{wp-1}$.  
\endproof

Note that the weight of a partition of $\frak{S}_{wp-1}$ must be less
than $w$, and all the decomposition numbers for blocks of weight $0,1$
or $2$ are known  [\textbf{JK}, \textbf{R}]. Therefore, Corollary
2.4 determines the decomposition numbers of $\frak{S}_{3p}$. Note,
too, that the proof fails when $w=4$ because $\lambda$ may have more
than one removable $r$-node in this case.  For example, suppose that
$p=3$ and consider the partition $\lambda =(6,4,1^{2})$, which has
the abacus configuration:

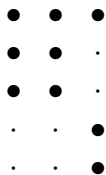
\[
\begin{array}{ccc}
\bigcirc & \bigcirc & \bigcirc \\ 
\bigcirc & \bigcirc & \cdot \\ 
\bigcirc & \bigcirc & \cdot \\ 
\cdot & \cdot & \bigcirc \\ 
\cdot & \cdot & \bigcirc
\end{array}
\]

\noindent 2.5 COROLLARY\ (Scopes [\textbf{S1}]) \textit{Suppose that
$B$ is a block of $\frak{S}_n$ such that the abacus configuration of every
partition in $B$ has the property that runner $i$ contains at least
$w$ more beads than runner $i-1$, for some $i$. Then each
decomposition number for $B$ can be equated with an explicit
decomposition number for some smaller symmetric group.}

\bigskip

\noindent\textit{Proof.} Making $w$ slides from the $p$-core, no position which we reach
allows us to move a bead on runner $i-1$ one space to the right. Therefore,
we may apply the Corollary 2.3, with $r=i$, to obtain the desired result.
\endproof

We remark that Scopes proved the stronger result that the block $B$ is
Morita equivalent to the block whose abacus configuration is obtained
by interchanging runners $i$ and $i-1$.

\[
\text{3. METHODS FOR ESTIMATING DECOMPOSITION NUMBERS} 
\]

\bigskip

We now present a collection of techniques for gathering information
about decomposition numbers. These ideas determine the decomposition numbers
for blocks of weight $0,1$ or $2$, and go some way in dealing with blocks of
higher weight. Many examples will appear later in this paper.

Suppose we are given a partition $\lambda$ and that we are trying to find 
$[S(\lambda):D(\mu)]$, for all $\mu$. We may assume that  $\lambda$
and $\mu$ are in the same block and that 
$\mu \trianglerighteq \lambda$, since otherwise $[S(\lambda):D(\mu)]=0$.  
 (Recall that $\mu\trianglerighteq\lambda$ if
$\sum_{i=1}^k\mu_i\ge\sum_{i=1}^k\lambda_i$, for all $k\ge1$. We say
that $\mu$ \textbf{dominates} $\lambda$.)
In particular, the number of (non-zero) parts of $\mu$ cannot
exceed the number of parts of $\lambda$. Hence, whatever abacus we
use to represent $\lambda$ we can  also use to represent $\mu$.
This follows because the number of parts of a partition can be read
off its abacus configuration by counting the number of beads after the
first gap.

 Also, recalling the definition of normal node from the last
section, observe that the normal $r$-nodes for $\mu$ can  also be
read off an abacus configuration for $\mu$ by considering the beads
on runners $r-1$ and $r$ in the abacus.

Our first rule is the abacus version of Corollary 2.2.

\bigskip

\noindent\textbf{Rule 1}.
Suppose that $\lambda$ has an abacus configuration such that
exactly $k$ beads on runner $r$ can be moved one space to the left and that 
none of the  beads on runner $r-1$ can be moved one space to the right, for some $r$. 
Then $[S(\lambda):D(\mu)]=[S(\bar{\lambda}):D(\bar{\mu})]$, where the abacus 
configuration for $\bar{\lambda}$ is obtained from that for $\lambda$ by moving 
to the left the $k$ possible beads on runner $r$, and the abacus
configuration for $\bar{\mu}$ is obtained by moving to the left the
$k$ beads on runner $r$ corresponding to the lowest $k$ normal nodes
in $\mu$.

\bigskip

\noindent\textbf{Notes}

(a) In practice, $\mu$ frequently has no addable $r$-nodes, so to
obtain $\bar{\mu}$ one simply moves to the left the $k$ possible beads
on runner $r$.

(b) If $k\ge1$ then Rule 1 equates $[S(\lambda):D(\mu)]$ with a
decomposition number of the same weight in a smaller symmetric group.

\bigskip

\noindent\textbf{Example}. If $\lambda$ and $\mu$ correspond to
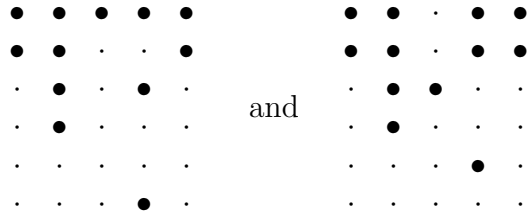
\[
\begin{array}{ccccc}
\bigcirc & \bigcirc & \bigcirc & \bigcirc & \bigcirc \\ 
\bigcirc & \bigcirc & \cdot & \cdot & \bigcirc \\ 
\cdot & \bigcirc & \cdot & \bigcirc & \cdot \\ 
\cdot & \bigcirc & \cdot & \cdot & \cdot \\ 
\cdot & \cdot & \cdot & \cdot & \cdot \\ 
\cdot & \cdot & \cdot & \bigcirc & \cdot
\end{array}
\text{ \ \ \ and\ \ \ } 
\begin{array}{ccccc}
\bigcirc & \bigcirc & \cdot & \bigcirc & \bigcirc \\ 
\bigcirc & \bigcirc & \cdot & \bigcirc & \bigcirc \\ 
\cdot & \bigcirc & \bigcirc & \cdot & \cdot \\ 
\cdot & \bigcirc & \cdot & \cdot & \cdot \\ 
\cdot & \cdot & \cdot & \bigcirc & \cdot \\ 
\cdot & \cdot & \cdot & \cdot & \cdot
\end{array}
\]
and $r=4$ then $\bar{\lambda}$ and $\bar{\mu}$ correspond to
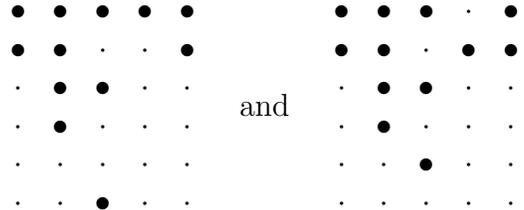
\[
\begin{array}{ccccc}
\bigcirc & \bigcirc & \bigcirc & \bigcirc & \bigcirc \\ 
\bigcirc & \bigcirc & \cdot & \cdot & \bigcirc \\ 
\cdot & \bigcirc & \bigcirc & \cdot & \cdot \\ 
\cdot & \bigcirc & \cdot & \cdot & \cdot \\ 
\cdot & \cdot & \cdot & \cdot & \cdot \\ 
\cdot & \cdot & \bigcirc & \cdot & \cdot
\end{array}
\text{\ \ \ and\ \ \ } 
\begin{array}{ccccc}
\bigcirc & \bigcirc & \bigcirc & \cdot & \bigcirc \\ 
\bigcirc & \bigcirc & \cdot & \bigcirc & \bigcirc \\ 
\cdot & \bigcirc & \bigcirc & \cdot & \cdot \\ 
\cdot & \bigcirc & \cdot & \cdot & \cdot \\ 
\cdot & \cdot & \bigcirc & \cdot & \cdot \\ 
\cdot & \cdot & \cdot & \cdot & \cdot
\end{array}
. 
\]
(Here, we could also apply Rule 1 with $r=2$, but not with $r=5$.)

\bigskip

\noindent\textbf{Rule 2}. Given a partition $\lambda$, 
Schaper's Theorem [\textbf{Sch}, \textbf{JM1}] gives us a linear
combination of Specht modules $S(\nu)$,  where 
$\nu\vartriangleright\lambda$ and $\nu$ belongs to the same block as
$\lambda$. If we know (for example, by induction) all of the
decomposition numbers for the Specht modules $S(\nu)$  appearing in
this sum then, in the Grothendieck group of $\frak S_n$, we can
rewrite this sum as a linear combination of irreducible modules
$D(\mu)$ with non-negative integer coefficients. Schaper's Theorem then
tells us that:

(a) if $D(\mu)$ appears in this linear combination with multiplicity $m>1$,
then $m\geq [S(\lambda):D(\mu)]\geq 1$; and

(b) if $D(\mu)$ appears in this linear combination with multiplicity 
$m\leq 1$, then $[S(\lambda):D(\mu)]=m$.

\bigskip

Note that Rule 2 gives us both upper and lower bounds on $[S(\lambda):D(\mu)]$.
Our next rule will provide another upper bound.

Suppose that $\mu$ has exactly $k$ normal $r$-nodes and let
$\bar{\mu}$  be the partition obtained from $\mu$ by removing these
nodes. Also, let $\Omega $ denote the set of partitions of $n-k$
which are obtained from $\lambda$ by removing $k$ $r$-nodes. Then
Kleshchev's Branching Theorem shows that
\[
[S(\lambda):D(\mu)]\leq \sum_{\omega \in \Omega }[S(\omega):D(\bar{\mu})]. 
\]
(Here, we interpret the right hand side to be zero when $\Omega $ is empty.)

\bigskip

\noindent\textbf{Rule 3} We may iterate the process just defined,
first removing all of the $k_{1}$ normal $r_{1}$-nodes from~$\mu$,
then taking all the $k_{2}$ normal $r_{2}$-nodes from the partition
$\bar\mu$, and so on, until we reach a stage where we can evaluate the
decomposition numbers on the right hand side of the inequality.

\bigskip

\noindent\textbf{Note} We do not increase the weight of the partitions
involved when we apply Rule 3. By this we mean that the weight of
$\bar{\mu}$ is at most the weight of $\mu$. To see this, first
observe that, in general, if the number of beads on runner $r-1$ is
$a$ and the number of beads on runner $r$ is $b,$ then moving a bead
left from runner $r$ to runner $r-1$ decreases the weight by $a-b+1$
(of course, a negative decrease corresponds to an increase). Hence, by
induction, moving $k $ beads left from runner $r$ to runner $r-1$
decreases the weight by $k(a-b+k).$ Now suppose that $\mu$ has
exactly $k$ normal $r$-nodes.  Then the definition of normal implies
that $k\geq b-a$; thus, $k(a-b+k)\geq 0$, so removing the $k$ normal
$r$-nodes does not increase the weight.

Observe that we can always apply Rule 3 to get an upper bound on
$[S(\lambda):D(\mu)]$ because at some point we will be able to
evaluate the right hand side of the inequality, if need be by
persevering until we reach the empty partition. If  in applying Rule~3
we remove $k_{1}$ normal $r_{1}$-nodes, $k_{2}$ normal $r_{2}$-nodes,
and so on, then we refer to $r_{1}^{k_{1}}r_{2}^{k_{2}}\dots $ as a
Kleshchev sequence for $\mu$.

\bigskip

The next Rule is due to the first author [\textbf{J1}]. 

\bigskip

\noindent\textbf{Rule 4} Assume that $\lambda$ and $\mu$ are partitions of $n$ with 
$\mu$ being $p$-regular, and $\lambda_{1}=\mu_{1}$. Let
$$
\bar\lambda=(\lambda_{2},\lambda_{3},\dots),\quad\text{and}
\qquad \bar\mu=(\mu_2,\mu_3,\dots).
$$
Then 
$[S(\lambda):D(\mu)]=[S(\bar\lambda):D(\bar\mu)]$.

\bigskip

This rule  says
that we can remove the first rows of $\lambda$ and $\mu$ without
changing the decomposition multiplicity $[S(\lambda):D(\mu)]$.
Analogously, we have the following Rule (see [\textbf{D},
\textbf{J1}]).

\bigskip

\noindent\textbf{Rule 5 }Assume that $\lambda$ and $\mu$ are partitions of $n$ with 
$\mu$ being $p$-regular, and that $\lambda$ and $\mu$ have the same first
column. Let

\[
\lambda^{(1)}=(\lambda_{1}-1,\lambda_{2}-1,\dots),\;\;\mu^{(1)}=(\mu
_{1}-1,\mu_{2}-1,\dots). 
\]
Then $[S(\lambda):D(\mu)]=[S(\lambda^{(1)}):D(\mu^{(1)})]$.

\bigskip

\noindent\textbf{Notes}

(a) Removing the first column from a partition corresponds to putting a bead
in the first gap.

(b) Removing the first row from a partition corresponds to removing the last
bead.

\bigskip

\noindent\textbf{Rule 6 }Assume that $\lambda$ is $p$-regular and we know the
decomposition numbers for every $S(\nu)$ with $\nu \trianglerighteq \lambda
$. Then we can express $D(\lambda)$ as a linear combination of the Specht
modules $S(\nu)$ with $\nu \trianglerighteq \lambda$. For all $r,$ the 
$r$-restriction of this linear combination of Specht modules is a module for 
$\frak{S}_{n-1}$.

\bigskip

We can apply Rule 6 to give an upper bound on a decomposition
number $[S(\lambda):D(\mu)]$ whenever we know the decomposition
numbers for every $S(\nu)$ with $\nu\vartriangleright\lambda$; see
the example at the end of section~5. Moreover, we can iterate this
process and perform a sequence $r_{1}^{k_{1}}r_{2}^{k_{2}}\dots $ of
restrictions, rather than just a single $r$-restriction.

Rule 6, which involves restricting an irreducible module, again gives us an
upper bound on decomposition numbers. We do not list the corresponding
result involving inducing an irreducible module, which would give us a lower
bound, for the following reason. \ If $[S(\lambda):D(\mu)]\geq 1$, then
inducing simple modules would \textit{perhaps }give this information, but
Rule 1 would \textit{certainly} give it.

\bigskip

\noindent\textbf{Rule 7 }Assume that $\lambda$ and $\mu$ are partitions of $n$ with 
$\mu$ being $p$-regular. Then $[S(\lambda):D(\mu)]=[S(\lambda^{\prime
}):D(\mu^{\ast })]$, where $\lambda^{\prime }$ is the conjugate of 
$\lambda$ and $\mu^{\ast }$ is the image of $\mu$ under the Mullineux map 
[\textbf{BO, FK}].

\bigskip

We reiterate that the Rules we have stated deal very well with many
decomposition numbers of blocks of small defect. Moreover, as we shall
see, Rules 1-6  add credibility to the conjecture which we discuss
next.

\bigskip

Our conjecture relates certain decomposition numbers for different
primes. 

Let $\lambda$ and $\mu$ be partitions, $\mu$ being $p$-regular, and
suppose that $\lambda$ and $\mu$ have the same $p$-core and the
same weight. Represent $\lambda$ on some abacus with $p$ runners. We
shall discuss the decomposition number $[S(\lambda):D(\mu)]$ so we
may assume that $\mu \trianglerighteq \lambda$; hence $\mu$ can be
represented on the same abacus as~$\lambda$. Suppose that $p^{\prime
}$ is a prime greater than $p$.  Our conjecture equates certain
$p^{\prime }$-modular decomposition numbers with $p$-modular
decomposition numbers. Let $\lambda^+$ denote the partition
obtained from $\lambda$ by adding $p^{\prime }-p$ empty runners to
the abacus (in any places) and let $\mu^+$ denote the partition
obtained from $\mu$ by adding $p^{\prime }-p$ empty runners to the
abacus configuration for $\mu$ (in the same places). We now put
forward the following Conjecture.

\bigskip

\noindent 3.1 CONJECTURE \textit{Suppose that }$p>w$\textit{. Then }
$[S(\lambda^+):D(\mu^+)]=[S(\lambda):D(\mu)]$.

\bigskip

Let $B$ be a block of $\mathfrak S_n$ of weight $w$. Then $B$ is a 
\textbf{block of small defect} if $p>w$.

Note that by Rule 1 the decomposition number
$[S(\lambda^+):D(\mu^+)]$ is independent of where the $p'-p$ empty
runners are inserted into the abacuses of $\lambda$ and $\mu$ (the
empty runners do, however, need to be in the same places).  Unless
otherwise stated we will assume that the abacuses for $\lambda^+$ and
$\mu^+$ are obtained from those for $\lambda$ and $\mu$, respectively,
by adding $p^{\prime }-p$ empty runners at the end.  

Rather than working with the symmetric group in characteristic~$p$ if,
instead, we work with the Hecke algebra of type $A$ at a complex $p$th
root of unity then our Conjecture is true, without any restriction on
$p$. This is part of the main result of our paper [\textbf{JM2}].

We remark that the assumption that $p>w$ in Conjecture~3.1 is
necessary. To see this let $\lambda=(3,1^2)$ and $\mu=(5)$ and take
$p=2$. Then $\lambda$ and $\mu$ are partitions of $2$--weight $2$ and
$[S(3,1^2):D(5)]=2$. These partitions have the following abacus
configurations.
$$\lambda=\begin{array}{cc} 
              \cdot & \bigcirc \\
              \bigcirc & \cdot \\
              \cdot & \bigcirc\\
              \cdot & \cdot \\
	   \end{array}
	   \hspace*{30mm}\mu=\begin{array}{cc}
              \bigcirc & \bigcirc \\
              \cdot & \cdot \\
              \cdot & \cdot \\
              \cdot & \bigcirc 
	    \end{array}$$
So we may take $\lambda^+=(5,2,1)$ and $\mu=(8)$ with $p'=3$ by
adding an empty right hand runner. However, $[S(5,2,1):D(8)]=1$
when $p'=3$. So $[S(\lambda):D(\mu)]\ne[S(\lambda^+):D(\mu^+)]$ in
this case.

% Although we believe that Conjecture 3.1 is true in the form stated
% above, there is slightly more evidence in favour of the following
% weaker conjecture (see Proposition 3.4).
% 
% \bigskip
% 
% $3.1'$ CONJECTURE \textit{Suppose that }$p>w$\textit{\ and some runner of the
% abacus contains no beads. Then }$[S(\lambda^+):D(\mu^+)]=[S(\lambda
% ):D(\mu)].$

\bigskip

We give further evidence in support of our conjecture after the
examples below.

\bigskip

\noindent\textbf{Example }Suppose that $p=5$ and $\lambda =(8,8,4,1)$ and 
$\mu =(12,9) $. Then we can represent $\lambda$ and~$\mu$ on an abacus
as follows.
\[
\lambda = 
\begin{array}{ccccc}
\cdot & \bigcirc & \cdot & \cdot & \cdot \\ 
\bigcirc & \cdot & \cdot & \cdot & \cdot \\ 
\bigcirc & \bigcirc & \cdot & \cdot & \cdot \\ 
\cdot & \cdot & \cdot & \cdot & \cdot \\ 
\cdot&\cdot&\cdot&\cdot& \cdot
\end{array}
\text{\quad and\quad }\mu = 
\begin{array}{ccccc}
\bigcirc & \bigcirc & \cdot & \cdot & \cdot \\ 
\cdot & \cdot & \cdot & \cdot & \cdot \\ 
\cdot & \bigcirc & \cdot & \cdot & \cdot \\ 
\bigcirc & \cdot & \cdot & \cdot & \cdot \\ 
\cdot&\cdot&\cdot&\cdot& \cdot
\end{array}
\]
Now let $p^{\prime }=7$. Then Rule 1 (applied $3$ times) ensures that 
$[S(\lambda^+):D(\mu^+)]$ is the same if
\[
\lambda^+= 
\begin{array}{ccccccc}
\cdot & \cdot & \cdot & \bigcirc & \cdot & \cdot & \cdot \\ 
\cdot & \bigcirc & \cdot & \cdot & \cdot & \cdot & \cdot \\ 
\cdot & \bigcirc & \cdot & \bigcirc & \cdot & \cdot & \cdot \\ 
\cdot & \cdot & \cdot & \cdot & \cdot & \cdot & \cdot \\ 
\cdot&\cdot&\cdot&\cdot&\cdot&\cdot& \cdot
\end{array}
\;\;\text{and\ \ }\mu^+= 
\begin{array}{ccccccc}
\cdot & \bigcirc & \cdot & \bigcirc & \cdot & \cdot & \cdot \\ 
\cdot & \cdot & \cdot & \cdot & \cdot & \cdot & \cdot \\ 
\cdot & \cdot & \cdot & \bigcirc & \cdot & \cdot & \cdot \\ 
\cdot & \bigcirc & \cdot & \cdot & \cdot & \cdot & \cdot \\ 
\cdot&\cdot&\cdot&\cdot&\cdot&\cdot& \cdot
\end{array}
\text{\ } 
\]
or if
\[
\lambda^+= 
\begin{array}{ccccccc}
\cdot & \bigcirc & \cdot & \cdot & \cdot & \cdot & \cdot \\ 
\bigcirc & \cdot & \cdot & \cdot & \cdot & \cdot & \cdot \\ 
\bigcirc & \bigcirc & \cdot & \cdot & \cdot & \cdot & \cdot \\ 
\cdot & \cdot & \cdot & \cdot & \cdot & \cdot & \cdot \\ 
\cdot&\cdot&\cdot&\cdot&\cdot&\cdot& \cdot
\end{array}
\;\;\text{and\ \ }\mu^+= 
\begin{array}{ccccccc}
\bigcirc & \bigcirc & \cdot & \cdot & \cdot & \cdot & \cdot \\ 
\cdot & \cdot & \cdot & \cdot & \cdot & \cdot & \cdot \\ 
\cdot & \bigcirc & \cdot & \cdot & \cdot & \cdot & \cdot \\ 
\bigcirc & \cdot & \cdot & \cdot & \cdot & \cdot & \cdot \\ 
\cdot&\cdot&\cdot&\cdot&\cdot&\cdot& \cdot
\end{array}
\]

Conjecture 3.1 says that $[S(\lambda^+):D(\mu^+)]=[S(\lambda):D(\mu)] $.  
Thus, in this example, our conjecture says that the decomposition multiplicity
\[
[S(2p-2,2p-2,p-1,1):D(3p-3,2p-1)] 
\]
is the same for $p=7$ as for $p=5$.  We know of no way to compute this multiplicity 
in general; however, using extensive computer calculations
L\"{u}beck and M\"{u}ller [\textbf{LM}, \textbf{LN}] have shown that if $p=5$ then
 $[S(2p-2,2p-2,p-1,1):D(3p-3,2p-1)]=1$.

If our conjecture is correct then it follows, as in the Example above, that
$$[S(2p-2,2p-2,p-1,1):D(3p-3,2p-1)]=1,\quad\text{for all\ } p>3.$$
On the other hand, if $[S(2p-2,2p-2,p-1,1):D(3p-3,2p-1)]\neq 1$, for
any $p>3$, then this provides a counterexample to the ``$pe>n$ Conjecture'' 
of [\textbf{J3}, Section 4].

We find it remarkable that if Conjecture 3.1 is correct then we can
produce a computer-free proof of L\"{u}beck and M\"{u}ller's result
 above. That is, we can deduce that $[S(8^{2},4,1):D(12,9)]=1$ when
$p=5$. Here is how this comes about. First, using Rules 1-7,
\[
1=[S(8,5,4,2):D(19)]\text{ when }p=5.
\]
Indeed, the decomposition matrices of $\frak{S}_{n}$ and $p=5 $ can be
calculated by hand for $n\leq 20$. Next,
$$\begin{array}{rcl@{\quad}l}
[S(8,5,4,2):D(19)]&=&[S(4^{2},3^{2},2,1^{3}):D(5^{3},4)],
                          &\text{when $p=5$, by Rule 7,}\\
                  &=&[S(8^{2},5^{2},4,1^{3}):D(9^{3},6)],
                          &\text{when $p=7$, if Conjecture 3.1 is correct.} \\
                  &=&[S(8,5^{3},4,2^{3}):D(11^{3})],
			  &\text{when $p=7$, by Rule 7.}
\end{array}
$$
Now, we are unable to evaluate the last decomposition number
when $p=7$ using Rules 1-7; however, using Rule 2 we can show that
$[S(8,5^{3},4,2^{3}):D(11^{3})]$ is either~$1$ or~$2$. By
investigating these two possibilities we can show that
 $[S(8,5^{3},4,2^{3}):D(11^{3})]=[S(10^{2},6,5,1^{2}):D(11^{3})]$,
when $p=7$, irrespective of the actual value of
$[S(8,5^{3},4,2^{3}):D(11^{3})]$. In turn,
$$\begin{array}{rcl@{\quad}l}
[S(10^{2},6,5,1^{2}):D(11^{3})]
 &=&[S(6^{2},4,3,1^{2}):D(7^{3})],&\text{when $p=5$, if Conjecture 3.1 is correct,}\\
 &=&[S(6,4^{2},3,2^{2}):D(12,9)], &\text{when $p=5$, by Rule 7.}
\end{array}$$
Again,  using Rules 1-7, we are only able to show that this last
decomposition number is either~$1$ or~$2$; however, by pursuing these
two possibilities in turn it can be shown that
\[
[S(6,4^{2},3,2^{2}):D(12,9)]=[S(8^{2},4,1):D(12,9)], \text{\ when $p=5$.} 
\]
Thus, by a very circuitous route, we have shown that if our conjecture is
true then one can deduce by hand that $[S(8^{2},4,1):D(12,9)]=1$
when $p=5$.

\bigskip

The argument for showing that $[S(8^{2},4,1):D(12,9)]=1$ when $p=5$,
consists of alternating applications of Conjecture 3.1 and Rule 7
(conjugation). In the absence of a proof of our conjecture, similar
arguments suggest that there are at least $p-3$ projective
indecomposable modules in blocks of weight $3$ which cannot, as yet,
be determined. These {are the blocks corresponding to the} $p$-cores
\[
(2^{p-2}),\,(3^{p-3}),\,(4^{p-4}),\,(5^{p-5}),\,\dots \,,(p-2)^{2}. 
\]
Notice that these cores occur in conjugate pairs, so there at least
$\frac{p-3}2$ independent decomposition numbers of weight 3 in characteristic
$p$ which current theory is unable to determine.

As evidence in support of our conjectures we present the following
Propositions which show that our conjectures are compatible with Rules
1--5. We remark that in the proofs of Propositions 3.2-3.6, the hypothesis that 
$p>w$ has immediate effect only in the proof of Proposition 3.3.

\bigskip

\noindent 3.2 PROPOSITION \textit{Assume that $p>w$ and that $\lambda$ and 
$\mu$ are partitions of $n$ of weight $w$ and with the same $p$-core and that
$\mu$ is $p$-regular. Assume that 
$[S(\alpha^+):D(\beta^+)]=[S(\alpha):D(\beta)]$ whenever~$\alpha$ and
$\beta$ are partitions of an integer less than $n$, with weight $w$.}

\textit{Suppose that $\lambda$ has an abacus configuration such that, for some
$r$, exactly $k>0$ beads on runner~$r$ can be moved one space to the
left and that none of the beads on runner $r-1$ can be moved one space
to the right, as in Rule 1.  Then
$[S(\lambda^+):D(\mu^+)]=[S(\lambda):D(\mu)]$.}

\bigskip

\noindent\textit{Proof.} Adopt the notation of Rule 1. Note that
$\bar{\lambda}$ and $\bar{\mu}$ have weight $w$. We may
assume that the abacuses for $\lambda^+$ and $\mu^+$ are obtained by
inserting $p'-p$ empty runners between runners $r$ and $r+1$ of the
abacuses for $\lambda$ and $\mu$, respectively; consequently,
 $(\bar{\lambda})^+=\overline{\lambda^+}$
 and $(\bar{\mu})^+=\overline{\mu^+}$. Therefore,
\begin{xalignat*}{2}
[ S(\lambda^+):D(\mu^+)]
   &=[S(\overline{\lambda^+}):D(\overline{\mu^+})],&&\text{by Rule 1,} \\
   &=[S((\bar{\lambda})^+):D((\bar{\mu})^+)] \\
   &=[S(\bar{\lambda}):D(\bar{\mu})],&&\text{by our induction hypothesis,} \\
   &=[S(\lambda):D(\mu)],&&\text{by Rule 1.}
\end{xalignat*}
\endproof

\bigskip

\noindent 3.3 PROPOSITION \textit{Assume that $p>w$ and that $\lambda$ and 
$\mu$ are partitions of $n$ of weight $w$ and with the same $p$-core and that
$\mu$ is $p$-regular. Assume that 
$[S(\nu^+):S(\mu^+)]=[S(\nu):D(\mu)]$
whenever $\nu\vartriangleright\lambda$ and  that by applying Rule~2 we
can deduce that $[S(\lambda):D(\mu)]\le m$. Then $[S(\lambda^+):D(\mu^+)]\le m$.
}

\bigskip

\noindent\textit{Proof.} Applying Schaper's Theorem to $S(\lambda)$ gives a
linear combination $\sum_\nu a_\nu S(\nu)$ of Specht modules $S(\nu)$,
where $a_\nu\ne0$ only if  $\lambda$ and $\nu$ belong to the same block
and $\nu\vartriangleright\lambda$.  Therefore, if $a_\nu\ne0$ then
$\lambda$ and $\nu$ can be represented on the same abacus and the
partitions $\nu$ that arise are determined by sliding beads up and
down the runners of an abacus for $\lambda$ in a specific way.
Moreover, the coefficient $a_\nu$ of a $S(\nu)$ in this linear
combination depends on the $p$-adic evaluation of the hook lengths
involved. Since $p>w$, no hook length in either of the partitions
$\lambda$ or $\lambda^+$ is divisible by $p^{2}$. Hence, Schaper's
Theorem applied to $S(\lambda^+)$ produces the linear combination
$\sum_\nu a_\nu S(\nu^+)$ of Specht modules. By  assumption, the
decomposition numbers $[S(\nu^+):D(\mu^+)]$ have already been
proved to be equal to $[S(\nu):D(\mu)]$. Therefore, the information
provided for $[S(\lambda):D(\mu)]$ by Rule 2 gives the same
information for $[S(\lambda^+):D(\mu^+)]$.
\endproof

\noindent 3.4 PROPOSITION \textit{Assume that $p>w$ and that $\lambda$
and $\mu$ are partitions of $n$ of weight $w$ and with the same
$p$-core and that $\mu$ is $p$-regular.  Suppose that
$[S(\alpha^+):D(\beta^+)]=[S(\alpha):D(\beta)]$ whenever $\alpha $ and
$\beta$ are partitions of an integer less than $n$, with weight at
most $w$. If Rule 3 gives 
$[S(\lambda):D(\mu)]\leq 
  \sum_{\omega \in \Omega}[S(\omega:D(\bar{\mu})]$ then 
$[S(\lambda^+):D(\mu^+)]
  \leq\sum_{\omega \in \Omega }[S(\omega):D(\bar{\mu})].$}

\bigskip

\noindent\textit{Proof.} Once again, we 
assume that the abacuses for $\lambda^+$ and $\mu^+$ are obtained by
inserting $p'-p$ empty runners between runners $r$ and $r+1$ of the
abacuses for $\lambda$ and $\mu$, respectively. Then 
$\overline{\mu^+}=(\bar{\mu})^+$ and hence
\begin{xalignat*}{2}
[ S(\lambda^+) :D(\mu^+)]&
  \leq \sum_{\omega \in \Omega}[S(\omega^+):D((\bar{\mu})^+)],&&\text{by Rule 3,} \\
  &=\sum_{\omega \in \Omega }[S(\omega):D(\bar{\mu})],
              &&\text{by our induction hypothesis.}
\end{xalignat*}
Note that we are justified in applying our induction hypothesis, in the
light of the note to Rule~3.
\endproof

\bigskip

\noindent 3.5 PROPOSITION \textit{Assume that $p>w$ and that $\lambda$
and $\mu$ are partitions of $n$ with the same
$p$-core and of weight $w,$ with $\mu$ $p$-regular. Assume, too, that
$\lambda$ and $\mu$ have the same first row, as in Rule 4.
Suppose that $[S(\alpha^+):D(\beta^+)]=[S(\alpha):D(\beta)]$
whenever $\alpha$ and $\beta$ are partitions of an integer
less than $n$, with weight at most $w$. Then 
$[S(\lambda^+):D(\mu^+)]=[S(\lambda):D(\mu)]$.}

\bigskip

\noindent\textit{Proof.} Since $\lambda$ and $\mu$ have the same first row, the
last bead on the abacus for $\lambda$ is in the same position as the last
bead for $\mu$. Removing this bead does not increase the weight.

Adopt the notation of Rule 4. Note that 
$(\lambda^{(1)})^+=(\lambda^+)^{(1)}$ and $(\mu^{(1)})^+=(\mu^+)^{(1)}$. 
Therefore,
\begin{xalignat*}{2}
[S(\lambda^+):D(\mu^+)]
  &=[S((\lambda^+)^{(1)}):D((\mu^+)^{(1)})], &&\text{by Rule 4,} \\
  &=[S((\lambda^{(1)})^+):D((\mu^{(1)})^+)] \\
  &=[S(\lambda^{(1)}):D(\mu^{(1)})],&&\text{by our induction hypothesis,} \\
  &=[S(\lambda):D(\mu)],&&\text{by Rule 4.}
\end{xalignat*}
\endproof

\bigskip
There is a more general version of row removal [\textbf{D}] which says
that if $\lambda_1+\dots+\lambda_s=\mu_1+\dots+\mu_s$, for some~$s$, then
$$[S(\lambda):D(\mu)]=[S(\lambda_1,\dots,\lambda_s):D(\mu_1,\dots,\mu_s)]
  [S(\lambda_{s+1},\lambda_{s+2},\dots):D(\mu_{s+1},\mu_{s+2},\dots)].$$
However, this result is not obviously compatible with Conjecture 3.1
when $s>1$ because it is easy to find examples where 
$\lambda^+_1+\dots+\lambda^+_s\ne\mu^+_1+\dots+\mu^+_s$.

The previous remark also applies for the general version of column
removal (Rule 5). Even so, we do have the following result.

\bigskip

\noindent 3.6 PROPOSITION \textit{Assume that $p>w$ and that $\lambda$
and $\mu$ are partitions of $n$ with the same
p-core and of weight $w$, with $\mu$ $p$-regular. Assume, too, that
$\lambda$ and $\mu$ have the same first column, as in Rule~5.
Suppose that $[S(\alpha^+):D(\beta^+)]=[S(\alpha):D(\beta)]$
whenever $\alpha$ and $\beta$ are partitions of an integer
less than $n$, with weight at most $w$.
Then $[S(\lambda^+):D(\mu^+)]=[S(\lambda):D(\mu)]$.}

\bigskip

\noindent\textit{Proof.} Since $\lambda$ and $\mu$ have the same first column, the
first gap in the abacus for $\lambda$ is in the same position as the first
gap in the abacus for $\mu$; say this is position $i.$ Suppose that
position $i$ on the $p$-abacus is position~$i^+$ on the 
$p^{\prime}$-abacus. Then, by repeated applications of Rule 5,
$[S(\lambda^+):D(\mu^+)]=[S(\alpha):D(\beta)]$, where $\alpha$
is obtained from $\lambda^+$ by filling the gaps up to and
including the gap at position $i^+$ and $\beta $ is obtained from
$\mu^+$ in the same way. Similarly,
$[S((\lambda^{(1)})^+):D((\mu^{(1)})^+)]=[S(\alpha):D(\beta)]$, 
where we adopt the notation of Rule 5. Therefore,
\begin{xalignat*}{2}
[S(\lambda^+):D(\mu^+)]
  &=[S((\lambda^{(1)})^+):D((\mu^{(1)})^+)] \\
  &=[S(\lambda^{(1)}):D(\mu^{(1)})],&&\text{ by our induction hypothesis,} \\
  &=[S(\lambda):D(\mu)],&&\text{ by Rule 5.}
\end{xalignat*}
\endproof

\bigskip

Roughly speaking, Propositions 3.2-3.6 say that if all decomposition numbers
were determined by Rules 1-6, then Conjecture 3.1 would be true by induction.

% \bigskip
% 
% \noindent3.7 COROLLARY \textit{Suppose that $p>2$ and that $\lambda$
% and $\mu$ are partitions of $n$ of $p$--weight $0$, $1$ or $2$. Then
% Conjecture 3.1 holds for $\lambda$ and $\mu$; that is,
% $[S(\lambda):D(\mu)]=[S(\lambda^+):D(\mu^+)]$.
% }
% 
% \bigskip
% 
% \textbf{CHECK THIS!!!!}
% 
% \noindent\textit{Proof.} When $\lambda$ and $\mu$ have weight $0$ or
% $1$ this is easily verified. When $\lambda$ and $\mu$ have $p$--weight $2$ then
% Scopes [\textbf{S2}] proved that $[S(\lambda):D(\mu)]\le1$ using
% Rules~1--5. Hence, $[S(\lambda^+):D(\mu^+)]\le1$ by Proposition~3.2.
% As the decomposition numbers are now known to be at most $1$ they are
% completely determined by Rule~2 (Schaper's theorem). The result
% follows.
% \endproof

\bigskip
Finally, we remark that Conjecture 3.1 is not obviously compatible
with Rule 6 because the restriction of a block of weight $w$ can have
arbitrarily large weight (in particular, the weight can be larger than
$p$). 
\bigskip

\[
\text{4. BLOCKS\ OF\ WEIGHT\ 3} 
\]

\bigskip

Now let $w=3$ and assume that $p>3$.  By repeatedly applying
Corollary 2.3 we can reduce the calculation of \textit{all}
decomposition numbers for blocks of weight 3 down to considering only
certain blocks or, equivalently, $p$-cores. We now describe the
abacuses for this minimal collection of $p$-cores. We assume,
without loss of generality, that  each of our abacuses has exactly
$3$ beads on runner $1$ and at least $3$ beads on every other runner.

\bigskip\noindent\textbf{CASE 1}

\noindent All of the runners contain exactly $3$ beads.

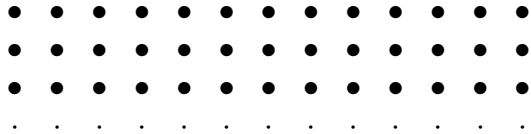
\[
\begin{array}{ccccccccccccc}
\bigcirc & \bigcirc & \bigcirc & \bigcirc & \bigcirc & \bigcirc & \bigcirc & 
\bigcirc & \bigcirc & \bigcirc & \bigcirc & \bigcirc & \bigcirc \\ 
\bigcirc & \bigcirc & \bigcirc & \bigcirc & \bigcirc & \bigcirc & \bigcirc & 
\bigcirc & \bigcirc & \bigcirc & \bigcirc & \bigcirc & \bigcirc \\ 
\bigcirc & \bigcirc & \bigcirc & \bigcirc & \bigcirc & \bigcirc & \bigcirc & 
\bigcirc & \bigcirc & \bigcirc & \bigcirc & \bigcirc & \bigcirc \\ 
\cdot&\cdot&\cdot&\cdot&\cdot&\cdot&\cdot&\cdot&\cdot&\cdot&\cdot&\cdot&\cdot\\ 
\end{array}
\]
Note that in this case the $p$-core is empty. In the light of Corollary 2.4
and the remark which follows it, we do not need to pursue Case 1 further.
(Case 1 is the only case which has to be considered if $w=1$.)

\bigskip\noindent\textbf{CASE 2}

\noindent The first $i-1$ runners contain exactly $3$ beads; runners $i$ up
to $j-1$ contain $4$ beads; and runners $j$ to $p$ contain $3$ beads.
Here,  $1<i<j\le p+1$, so there are $\binom p2$ such $p$-cores.

\bigskip 
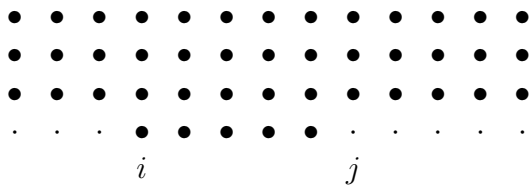
\[
\begin{array}{ccccccccccccc}
\bigcirc & \bigcirc & \bigcirc & \bigcirc & \bigcirc & \bigcirc & \bigcirc & 
\bigcirc & \bigcirc & \bigcirc & \bigcirc & \bigcirc & \bigcirc \\ 
\bigcirc & \bigcirc & \bigcirc & \bigcirc & \bigcirc & \bigcirc & \bigcirc & 
\bigcirc & \bigcirc & \bigcirc & \bigcirc & \bigcirc & \bigcirc \\ 
\bigcirc & \bigcirc & \bigcirc & \bigcirc & \bigcirc & \bigcirc & \bigcirc & 
\bigcirc & \bigcirc & \bigcirc & \bigcirc & \bigcirc & \bigcirc \\ 
\cdot&\cdot&\cdot& \bigcirc & \bigcirc & \bigcirc & \bigcirc & \bigcirc &\cdot&\cdot&\cdot&\cdot& \cdot\\ 
 & & & i & & & & & j & & & & \\ 
\end{array}
\]
Note that in this case the $p$-core is $(i-1)^{j-i},$ a partition of 
$ij-i^{2}+i-j$. (Cases 1 and 2 are the only cases which have to be considered
if $w=2$.)

\bigskip\noindent\textbf{CASE 3}

\noindent The first $i-1$ runners contain exactly $3$ beads; runner
$i$ contains $4$ beads; runners $i+1$ to $j-1$ contain $5$ beads;
runners $j$ to $k-1$ contain $4$ beads; and runners $k$ to $p$
contain $3$ beads. Here, we allow $2<i+1<j\le k\le p+1$, so there are 
$\binom{p-1}3+\binom{p-1}2=\binom p3$ such $p$-cores.

\bigskip 
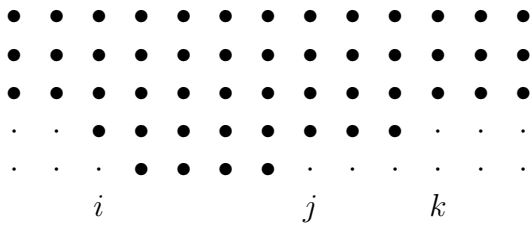
\[
\begin{array}{ccccccccccccc}
\bigcirc & \bigcirc & \bigcirc & \bigcirc & \bigcirc & \bigcirc & \bigcirc & 
\bigcirc & \bigcirc & \bigcirc & \bigcirc & \bigcirc & \bigcirc \\ 
\bigcirc & \bigcirc & \bigcirc & \bigcirc & \bigcirc & \bigcirc & \bigcirc & 
\bigcirc & \bigcirc & \bigcirc & \bigcirc & \bigcirc & \bigcirc \\ 
\bigcirc & \bigcirc & \bigcirc & \bigcirc & \bigcirc & \bigcirc & \bigcirc & 
\bigcirc & \bigcirc & \bigcirc & \bigcirc & \bigcirc & \bigcirc \\ 
\cdot&\cdot& \bigcirc & \bigcirc & \bigcirc & \bigcirc & \bigcirc & \bigcirc & 
\bigcirc & \bigcirc &\cdot&\cdot&\cdot\\ 
\cdot&\cdot&\cdot& \bigcirc & \bigcirc & \bigcirc & \bigcirc &\cdot&\cdot&\cdot&\cdot&\cdot&\cdot\\ 
 & & i & & & & & j & & & k & &
\end{array}
\]
Note that in this case the $p$-core is $((p-k+2i)^{j-i-1},(i-1)^{k-i})$.

\bigskip\noindent\textbf{CASE 4}

\noindent The first $i-1$ runners contain exactly $3$ beads; runners
$i$ to $j-1$ contain $5$ beads; runners $j$ to $k-1$ contain $4$
beads; and runners $k$ to $p$ contain $3$ beads. Here, we allow
$1<i<j\le k\le p+1$, so there are 
$\binom p3+\binom p2=\binom{p+1}3$ such $p$-cores.  

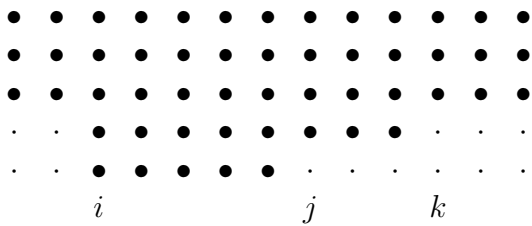
\[
\begin{array}{ccccccccccccc}
\bigcirc & \bigcirc & \bigcirc & \bigcirc & \bigcirc & \bigcirc & \bigcirc & 
\bigcirc & \bigcirc & \bigcirc & \bigcirc & \bigcirc & \bigcirc \\ 
\bigcirc & \bigcirc & \bigcirc & \bigcirc & \bigcirc & \bigcirc & \bigcirc & 
\bigcirc & \bigcirc & \bigcirc & \bigcirc & \bigcirc & \bigcirc \\ 
\bigcirc & \bigcirc & \bigcirc & \bigcirc & \bigcirc & \bigcirc & \bigcirc & 
\bigcirc & \bigcirc & \bigcirc & \bigcirc & \bigcirc & \bigcirc \\ 
\cdot&\cdot& \bigcirc & \bigcirc & \bigcirc & \bigcirc & \bigcirc & \bigcirc & 
\bigcirc & \bigcirc &\cdot&\cdot&\cdot\\ 
\cdot&\cdot& \bigcirc & \bigcirc & \bigcirc & \bigcirc & \bigcirc &\cdot&\cdot&\cdot&\cdot&\cdot& 
\cdot\\ 
 & & i & & & & & j & & & k & &
\end{array}
\]
Note that in this case the $p$-core is $((p-k+2i-1)^{j-i},(i-1)^{k-i})$.
 (The reader should have no difficulty working out the additional cases
which arise for $w=4$, and for higher weights.)

\bigskip
 Thus, in general there are 
$\binom{p+1}3+\binom p3+\binom p2+1=2\binom{p+1}3+1$ different
$p$-cores which need to be considered. In each case, the decomposition
numbers for the different cores are often very similar (depending on
the parameters $i,j,k$); however, the explosion of delicate subcases
makes it very difficult to write down a convincing argument for
general $p$.

\bigskip

Of the four cases that need to be considered when $w=3$, Case 3
was overlooked in [\textbf{MR}], thereby further jeopardizing their
claim that the decomposition numbers for $w=3$ can be determined and
that they all have value $0$ or $1$. We now apply our
methods to obtain certain decomposition numbers for $w=3$ in Cases 2,
3 and 4. Some of these results already appear in
[\textbf{MR}]. 

Suppose that we have fixed a $p$-core $\rho$ \textit{and} an
abacus configuration for $\rho$ as above. Then, as in
[\textbf{MR}, \textbf{S2}], we use the following notation  for the
partitions of weight $3$ with $p$-core $\rho$.
\begin{enumerate}
\item Let $\langle r\rangle $ be the partition  whose
    abacus is obtained by moving the last bead on runner $r$
    (of the abacus for $\rho$), down $3$ places.

\item Let $\langle r^{2}\rangle $ be the partition
    whose abacus is obtained by moving the last bead on runner
    $r$ down $2$ places and the second last bead on runner $r$,
    down one place.

\item Let $\langle r^{3}\rangle $ be the partition
    whose abacus is obtained by moving the last $3$ beads on runner $r$,
    down one place each.

\item For $r\neq s$, let $\langle r,s\rangle $ be the
    partition whose abacus is obtained by moving the last bead on runner
    $r$ down $2$ places and the last bead on runner $s$ down one place.

\item For $r\neq s$, let $\langle r^{2},s\rangle $ be
    the partition whose abacus is obtained by moving the last $2$ beads on
    runner $r$ down one place each, and the last bead on runner $s$ down one
    place.

\item For $r,s,t$ distinct, let $\langle r,s,t\rangle $
    be the partition whose abacus is obtained by moving 
    the last bead on runners $r$, $s$ and $t$ down one place each.
\end{enumerate}

\bigskip

\noindent\textbf{Some decomposition numbers in Case 2}

\noindent Assume that the $p$-core belongs to Case 2. That is,
runners up to runner $i-1$ contain $3$ beads; runner $i$ contains
$4$ beads; after this, there are some or no runners with $4$ beads;
any remaining runners contain $3$ beads.  Let runner $j$ be the first
runner after runner $i$ with $3$ beads.

By Rule 1, we can equate the decomposition number 
$[S(\lambda):D(\mu)]$ with a decomposition number of weight $3$ in a
smaller symmetric group for all partitions $\lambda$ and $\mu$
in this block except for  when $\lambda$ is one of the partitions
$\alpha^{\ast}$,     $\beta^{\ast }$,    $\gamma^{\ast }$,         
$\alpha^{\natural}$, $\beta^{\natural}$, $\gamma^{\natural }$, 
$\alpha^{(u)}$,      $\beta^{(u)}$  or   $\gamma^{(u)}$ 
where
\[
\alpha^{\ast }=\langle i^{2}\rangle ,\quad
\beta^{\ast }=\langle i,i-1\rangle ,\quad
\gamma^{\ast }=\langle i-1\rangle 
\]
and
\[
\alpha^{\natural }=\langle i^{3}\rangle ,\quad
\beta^{\natural }=\langle (i-1)^{2},i\rangle ,\quad
\gamma^{\natural }=\langle (i-1)^{2}\rangle 
\]
and
\[
\alpha^{(u)}=\langle i^{2},u\rangle,  \quad \text{ for }1\leq u\leq p
\text{ and }u\neq i-i,i 
\]
(Rule 1 deals with the case $u=i-1,$ and $u=i$ gives $\alpha^{\natural }$),
\[
\beta^{(u)}=\langle i-1,i,u\rangle, \quad\text{for }1\leq u\leq p
\text{ and }u\neq i-1,i 
\]
(Rule 1 deals with the case $u=i-1,$ and $u=i$ gives $\beta^{\ast }$),
and
\[
\gamma^{(u)}=\langle i-1,u\rangle ,\quad\text{for }1\leq u\leq p
\text{ and }u\neq i-1,i. 
\]
(Rule 1 deals with the case $u=i,$ and $u=i-1$ gives $\gamma^{\ast }$).

We call 
$\alpha^{\ast },\beta^{\ast },\gamma^{\ast },
 \alpha^{\natural },\beta^{\natural },\gamma^{\natural }$,
 $\alpha^{(u)},\beta^{(u)},\gamma^{(u)}$ 
the \textit{exceptional partitions} for Case 2.  The abacuses for the
exceptional partitions in Case 2 are as follows:
\[
\alpha^{\ast }= 
\begin{array}{lllllllllllll}
\bigcirc & \bigcirc & \bigcirc & \bigcirc & \bigcirc & \bigcirc & \bigcirc & 
\bigcirc & \bigcirc & \bigcirc & \bigcirc & \bigcirc & \bigcirc \\ 
\bigcirc & \bigcirc & \bigcirc & \bigcirc & \bigcirc & \bigcirc & \bigcirc & 
\bigcirc & \bigcirc & \bigcirc & \bigcirc & \bigcirc & \bigcirc \\ 
\bigcirc & \bigcirc & \bigcirc &\cdot& \bigcirc & \bigcirc & \bigcirc & \bigcirc
     & \bigcirc & \bigcirc & \bigcirc & \bigcirc & \bigcirc \\ 
\cdot&\cdot&\cdot& \bigcirc & \bigcirc & \bigcirc & \bigcirc & \bigcirc &\cdot&\cdot&\cdot&\cdot& 
\cdot\\ 
\cdot&\cdot&\cdot&\cdot&\cdot&\cdot&\cdot&\cdot&\cdot&\cdot&\cdot&\cdot&\cdot\\ 
\cdot&\cdot&\cdot& \bigcirc &\cdot&\cdot&\cdot&\cdot&\cdot&\cdot&\cdot&\cdot& \cdot
\end{array}
\]
so $\alpha^{\ast }=(2p-j+2i,i^{j-i},1^{p-i})$;
\[
\beta^{\ast }=
\begin{array}{lllllllllllll}
\bigcirc  & \bigcirc  & \bigcirc  & \bigcirc  & \bigcirc  & \bigcirc  & 
\bigcirc  & \bigcirc  & \bigcirc  & \bigcirc  & \bigcirc  & \bigcirc  & 
\bigcirc  \\ 
\bigcirc  & \bigcirc  & \bigcirc  & \bigcirc  & \bigcirc  & \bigcirc  & 
\bigcirc  & \bigcirc  & \bigcirc  & \bigcirc  & \bigcirc  & \bigcirc  & 
\bigcirc  \\ 
\bigcirc  & \bigcirc  &\cdot& \bigcirc  & \bigcirc  & \bigcirc  & \bigcirc  & 
\bigcirc  & \bigcirc  & \bigcirc  & \bigcirc  & \bigcirc  & \bigcirc  \\ 
\cdot&\cdot& \bigcirc  &\cdot& \bigcirc  & \bigcirc  & \bigcirc  & \bigcirc  &\cdot&\cdot&\cdot& 
\cdot&\cdot\\ 
\cdot&\cdot&\cdot&\cdot&\cdot&\cdot&\cdot&\cdot&\cdot&\cdot&\cdot&\cdot&\cdot\\ 
\cdot&\cdot&\cdot& \bigcirc  &\cdot&\cdot&\cdot&\cdot&\cdot&\cdot&\cdot&\cdot& \cdot
\end{array}
\]
so $\beta^{\ast }=(2p-j+2i,i^{j-i-1},i-1,1^{p-i+1})$;
\[
\gamma^{\ast }= 
\begin{array}{lllllllllllll}
\bigcirc & \bigcirc & \bigcirc & \bigcirc & \bigcirc & \bigcirc & \bigcirc & 
\bigcirc & \bigcirc & \bigcirc & \bigcirc & \bigcirc & \bigcirc \\ 
\bigcirc & \bigcirc & \bigcirc & \bigcirc & \bigcirc & \bigcirc & \bigcirc & 
\bigcirc & \bigcirc & \bigcirc & \bigcirc & \bigcirc & \bigcirc \\ 
\bigcirc & \bigcirc &\cdot& \bigcirc & \bigcirc & \bigcirc & \bigcirc & \bigcirc
     & \bigcirc & \bigcirc & \bigcirc & \bigcirc & \bigcirc \\ 
\cdot&\cdot&\cdot& \bigcirc & \bigcirc & \bigcirc & \bigcirc & \bigcirc &\cdot&\cdot&\cdot&\cdot& 
\cdot\\ 
\cdot&\cdot&\cdot&\cdot&\cdot&\cdot&\cdot&\cdot&\cdot&\cdot&\cdot&\cdot&\cdot\\ 
\cdot&\cdot& \bigcirc &\cdot&\cdot&\cdot&\cdot&\cdot&\cdot&\cdot&\cdot&\cdot& \cdot
\end{array}
\]
so $\gamma^{\ast }=(2p-j+2i-1,i^{j-i},1^{p-i+1})$;
\[
\alpha^{\natural }= 
\begin{array}{lllllllllllll}
\bigcirc & \bigcirc & \bigcirc & \bigcirc & \bigcirc & \bigcirc & \bigcirc & 
\bigcirc & \bigcirc & \bigcirc & \bigcirc & \bigcirc & \bigcirc \\ 
\bigcirc & \bigcirc & \bigcirc &\cdot& \bigcirc & \bigcirc & \bigcirc & \bigcirc
     & \bigcirc & \bigcirc & \bigcirc & \bigcirc & \bigcirc \\ 
\bigcirc & \bigcirc & \bigcirc & \bigcirc & \bigcirc & \bigcirc & \bigcirc & 
\bigcirc & \bigcirc & \bigcirc & \bigcirc & \bigcirc & \bigcirc \\ 
\cdot&\cdot&\cdot& \bigcirc & \bigcirc & \bigcirc & \bigcirc & \bigcirc &\cdot&\cdot&\cdot&\cdot& 
\cdot\\ 
\cdot&\cdot&\cdot& \bigcirc &\cdot&\cdot&\cdot&\cdot&\cdot&\cdot&\cdot&\cdot&\cdot\\ 
\cdot&\cdot&\cdot&\cdot&\cdot&\cdot&\cdot&\cdot&\cdot&\cdot&\cdot&\cdot& \cdot
\end{array}
\]
so $\alpha^{\natural }=(p-j+2i,i^{j-i},1^{2p-i})$;
\[
\beta^{\natural }= 
\begin{array}{lllllllllllll}
\bigcirc & \bigcirc & \bigcirc & \bigcirc & \bigcirc & \bigcirc & \bigcirc & 
\bigcirc & \bigcirc & \bigcirc & \bigcirc & \bigcirc & \bigcirc \\ 
\bigcirc & \bigcirc &\cdot& \bigcirc & \bigcirc & \bigcirc & \bigcirc & \bigcirc
     & \bigcirc & \bigcirc & \bigcirc & \bigcirc & \bigcirc \\ 
\bigcirc & \bigcirc & \bigcirc & \bigcirc & \bigcirc & \bigcirc & \bigcirc & 
\bigcirc & \bigcirc & \bigcirc & \bigcirc & \bigcirc & \bigcirc \\ 
\cdot&\cdot& \bigcirc &\cdot& \bigcirc & \bigcirc & \bigcirc & \bigcirc &\cdot&\cdot&\cdot&\cdot& 
\cdot\\ 
\cdot&\cdot&\cdot& \bigcirc &\cdot&\cdot&\cdot&\cdot&\cdot&\cdot&\cdot&\cdot&\cdot\\ 
\cdot&\cdot&\cdot&\cdot&\cdot&\cdot&\cdot&\cdot&\cdot&\cdot&\cdot&\cdot& \cdot
\end{array}
\]
so $\beta^{\natural }=(p-j+2i,i^{j-i-1},i-1,1^{2p-i+1})$;
\[
\gamma^{\natural }= 
\begin{array}{lllllllllllll}
\bigcirc & \bigcirc & \bigcirc & \bigcirc & \bigcirc & \bigcirc & \bigcirc & 
\bigcirc & \bigcirc & \bigcirc & \bigcirc & \bigcirc & \bigcirc \\ 
\bigcirc & \bigcirc &\cdot& \bigcirc & \bigcirc & \bigcirc & \bigcirc & \bigcirc
     & \bigcirc & \bigcirc & \bigcirc & \bigcirc & \bigcirc \\ 
\bigcirc & \bigcirc & \bigcirc & \bigcirc & \bigcirc & \bigcirc & \bigcirc & 
\bigcirc & \bigcirc & \bigcirc & \bigcirc & \bigcirc & \bigcirc \\ 
\cdot&\cdot&\cdot& \bigcirc & \bigcirc & \bigcirc & \bigcirc & \bigcirc &\cdot&\cdot&\cdot&\cdot& 
\cdot\\ 
\cdot&\cdot& \bigcirc &\cdot&\cdot&\cdot&\cdot&\cdot&\cdot&\cdot&\cdot&\cdot&\cdot\\ 
\cdot&\cdot&\cdot&\cdot&\cdot&\cdot&\cdot&\cdot&\cdot&\cdot&\cdot&\cdot& \cdot
\end{array}
\]
so $\gamma^{\natural }=(p-j+2i-1,i^{j-i},1^{2p-i+1})$;
\begin{eqnarray*}
\alpha^{(u)} &=& 
\begin{array}{lllllllllllll}
\bigcirc & \bigcirc & \bigcirc & \bigcirc & \bigcirc & \bigcirc & \bigcirc & 
\bigcirc & \bigcirc & \bigcirc & \bigcirc & \bigcirc & \bigcirc \\ 
\bigcirc & \bigcirc & \bigcirc & \bigcirc & \bigcirc & \bigcirc & \bigcirc & 
\bigcirc & \bigcirc & \bigcirc & \bigcirc & \bigcirc & \bigcirc \\ 
\bigcirc & \bigcirc & \bigcirc &\cdot& \bigcirc & \bigcirc & \bigcirc & \bigcirc
     & \bigcirc & \bigcirc & \bigcirc & \bigcirc & \bigcirc \\ 
\cdot&\cdot&\cdot& \bigcirc & \bigcirc & \bigcirc & \bigcirc & \bigcirc &\cdot&\cdot&\cdot&\cdot& 
\cdot\\ 
\cdot&\cdot&\cdot& \bigcirc &\cdot&\cdot&\cdot&\cdot&\cdot&\cdot&\cdot&\cdot&\cdot\\ 
\cdot&\cdot&\cdot&\cdot&\cdot&\cdot&\cdot&\cdot&\cdot&\cdot&\cdot&\cdot& \cdot
\end{array}
\\
&&\qquad+\text{ one move on runner }u,
\end{eqnarray*}
so 
$\alpha^{(u)}=\begin{cases}
    (p-j+2i,u-j+i+1,(i+1)^{j-i},2^{p-u},1^{u-i-1}),   &\text{if\ } j\leq u\leq p,\\
    (p-j+i+u,p-j+2i+1,(i+1)^{j-u-1},i^{u-i},1^{p-1}), &\text{if\ } i<u<j,\\
    (p-j+2i,i^{j-i},u+1,2^{p-i},1^{i-u-1}),           &\text{if\ } 1\leq u<i-1;
\end{cases}$
\begin{eqnarray*}
\beta^{(u)} &=& 
\begin{array}{lllllllllllll}
\bigcirc & \bigcirc & \bigcirc & \bigcirc & \bigcirc & \bigcirc & \bigcirc & 
\bigcirc & \bigcirc & \bigcirc & \bigcirc & \bigcirc & \bigcirc \\ 
\bigcirc & \bigcirc & \bigcirc & \bigcirc & \bigcirc & \bigcirc & \bigcirc & 
\bigcirc & \bigcirc & \bigcirc & \bigcirc & \bigcirc & \bigcirc \\ 
\bigcirc & \bigcirc &\cdot& \bigcirc & \bigcirc & \bigcirc & \bigcirc & \bigcirc
     & \bigcirc & \bigcirc & \bigcirc & \bigcirc & \bigcirc \\ 
\cdot&\cdot& \bigcirc &\cdot& \bigcirc & \bigcirc & \bigcirc & \bigcirc &\cdot&\cdot&\cdot&\cdot& 
\cdot\\ 
\cdot&\cdot&\cdot& \bigcirc &\cdot&\cdot&\cdot&\cdot&\cdot&\cdot&\cdot&\cdot&\cdot\\ 
\cdot&\cdot&\cdot&\cdot&\cdot&\cdot&\cdot&\cdot&\cdot&\cdot&\cdot&\cdot& \cdot
\end{array}
\\
&&\qquad+\text{ one move on runner }u,
\end{eqnarray*}
so 
$\beta^{(u)}=
\begin{cases}
    (p-j+2i,u-j+i+1,(i+1)^{j-i-1},i,2^{p-u},1^{u-i}),    &\text{if\ } j\leq u\leq p,\\
    (p-j+i+u,p-j+2i+1,(i+1)^{j-u-1},i^{u-i-1},i-1,1^{p-i+1}),&\text{if\ } i<u<j,\\
    (p-j+2i,i^{j-i-1},i-1,u+1,2^{p-i+1},1^{i-u-2}),          &\text{if\ }1\leq u<i-1;
\end{cases}$
\begin{eqnarray*}
\gamma^{(u)} &=& 
\begin{array}{lllllllllllll}
\bigcirc & \bigcirc & \bigcirc & \bigcirc & \bigcirc & \bigcirc & \bigcirc & 
\bigcirc & \bigcirc & \bigcirc & \bigcirc & \bigcirc & \bigcirc \\ 
\bigcirc & \bigcirc & \bigcirc & \bigcirc & \bigcirc & \bigcirc & \bigcirc & 
\bigcirc & \bigcirc & \bigcirc & \bigcirc & \bigcirc & \bigcirc \\ 
\bigcirc & \bigcirc &\cdot& \bigcirc & \bigcirc & \bigcirc & \bigcirc & \bigcirc
     & \bigcirc & \bigcirc & \bigcirc & \bigcirc & \bigcirc \\ 
\cdot&\cdot&\cdot& \bigcirc & \bigcirc & \bigcirc & \bigcirc & \bigcirc &\cdot&\cdot&\cdot&\cdot& 
\cdot\\ 
\cdot&\cdot& \bigcirc &\cdot&\cdot&\cdot&\cdot&\cdot&\cdot&\cdot&\cdot&\cdot&\cdot\\ 
\cdot&\cdot&\cdot&\cdot&\cdot&\cdot&\cdot&\cdot&\cdot&\cdot&\cdot&\cdot& \cdot
\end{array}
\\
&&\qquad+\text{ one move on runner }u,
\end{eqnarray*}
so 
$\gamma^{(u)}=
\begin{cases}
    (p-j+2i-1,u-j+i+1,(i+1)^{j-i},2^{p-u},1^{u-i}),  &\text{if\ }j\leq u\leq p,\\
    (p-j+u+i,p-j+2i,(i+1)^{j-u-1},i^{u-i},1^{p-i+1}),&\text{if\ }i<u<j,\\
    (p-j+2i-1,i^{j-i},u+1,2^{p-i+1},1^{i-u-2}),      &\text{if\ }1\leq u<i-1.
\end{cases}$

\bigskip

\noindent\textbf{Note}
The  three partitions $\alpha^{\natural },\beta^{\natural },\gamma^{\natural }$
are $p$-singular and, moreover, $\alpha^{\ast }\vartriangleright \beta^{\ast
}\vartriangleright \gamma^{\ast }$ and $\alpha^{\natural
}\vartriangleright \beta^{\natural }\vartriangleright \gamma^{\natural }$
and $\alpha^{(u)}\vartriangleright \beta^{(u)}\vartriangleright \gamma
^{(u)}$ for all $u$ with $1\leq u\leq p$ and $u\neq i-1,i$.
Also, if $\lambda \in \{\alpha ,\beta ,\gamma \}$ then
$$\lambda^{\ast}\vartriangleright\lambda^{(j-1)}
                 \vartriangleright\lambda^{(j-2)}
                 \vartriangleright\dots\vartriangleright \lambda^{(i+1)}
                 \vartriangleright\lambda^{(p)}
		 \vartriangleright\lambda^{(p-1)}
                 \vartriangleright\dots\vartriangleright \lambda^{(j)}
                 \vartriangleright\lambda^{(i-2)}
		 \vartriangleright\lambda^{(i-3)}\dots\vartriangleright \lambda^{(1)}
                 \vartriangleright\lambda^{\natural }.$$

\bigskip

\noindent 4.1 PROPOSITION \textit{Assume that $1\leq v\leq p$ and $v\neq i-1,i$.
\begin{enumerate}
    \item For $\lambda$\ arbitrary and 
	$\mu\in\{\alpha^{\ast},\alpha^{(v)}\}$ we can compute 
	$[S(\lambda):D(\mu)]$.
   \item For $\lambda$ an exceptional partition and 
       $\mu\in\{\beta^{\ast},\beta^{(v)}\}$ we can compute 
       $[S(\lambda):D(\mu)]$.
   \item For $\lambda$ an exceptional partition and 
       $\mu\in\{\gamma^{\ast},\gamma^{(v)}\}$ we can compute 
       $[S(\lambda):D(\mu)]$.
\end{enumerate}
}

\noindent\textit{Proof.} Recall that the decomposition numbers are
known for blocks of weight $0,1$ and $2$. We prove that we can reduce
the calculation of the decomposition numbers in the Proposition to one
of these cases.

(1) Suppose that 
$\mu \in \{\alpha^{\ast },\alpha^{(v)}\}$. Then $\mu$ has $2$
normal $i$-nodes and every $\lambda$ in the same block as $\mu$ has
at most $2$ removable $i$-nodes. Therefore, we can apply Proposition 2.1 to 
compute $[S(\lambda):D(\mu)]$.

(2) Suppose that $\mu \in \{\beta^{\ast },\beta^{(v)}\}$.

If $i\neq 2$ and $\mu \neq \beta^{(i-2)}$ then $\mu$ has exactly one
normal $(i-1)$-node, and $\lambda$ has at most $1$ removable $(i-1)$-node,
so we can apply Proposition 2.1 again.

Assume that $\mu =\beta^{(i-2)}$ and $i\neq 2,3$. Then $\mu$ has exactly
one normal $(i-2)$-node, and $\lambda$ has at most $1$ removable 
$(i-2)$-node, so we can apply Proposition 2.1 again. Note that if 
$\mu =\beta^{(i-2)}$ and $i=3,$ then $\mu$ is $p$-singular.

Assume that $i=2$ and $\mu$ is $p$-regular. Then $\mu =\beta^{(v)}$
for some $v$ with $j\leq v\leq p$. We need only consider those
partitions $\lambda$ for which the first part of $\mu$ 
is larger than the first part of $\lambda$ (since, otherwise, 
either $\mu\ntrianglerighteq \lambda$ or we can apply row removal).
Therefore, 
$\lambda \in \{\gamma^{\natural }, \gamma^{(j)},\dots,\gamma^{(p)}\}$. 
But the first columns of $\beta^{(v)}$ and $\gamma^{(u)}$ have the
same length and we can apply Rule~5. Also, unless $v=j$, we see
that $\beta^{(v)}$ has a normal $v$-node while $\gamma^{\natural }$
has no removable $v$-node; so, $[S(\lambda):D(\mu)]=0$ by Proposition 2.1.

We are now left with one final case, namely, $i=2,\;\mu =\beta^{(j)}$
(with $j\leq p$, since otherwise $\mu$ is $p$-singular) and
$\lambda =\gamma^{\natural }$.

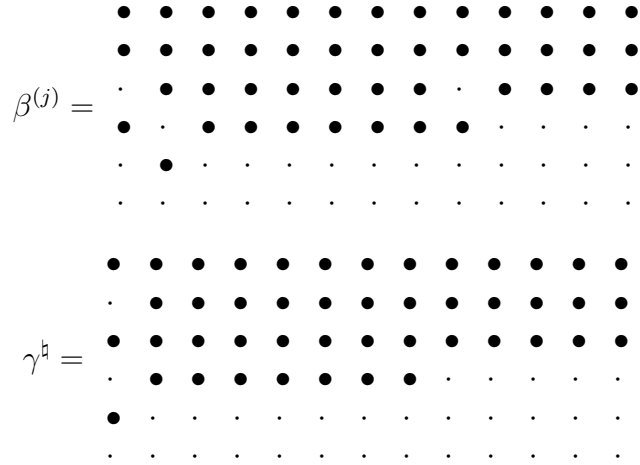
\[
\beta^{(j)}= 
\begin{array}{lllllllllllll}
\bigcirc & \bigcirc & \bigcirc & \bigcirc & \bigcirc & \bigcirc & \bigcirc & 
\bigcirc & \bigcirc & \bigcirc & \bigcirc & \bigcirc & \bigcirc \\ 
\bigcirc & \bigcirc & \bigcirc & \bigcirc & \bigcirc & \bigcirc & \bigcirc & 
\bigcirc & \bigcirc & \bigcirc & \bigcirc & \bigcirc & \bigcirc \\ 
\cdot& \bigcirc & \bigcirc & \bigcirc & \bigcirc & \bigcirc & \bigcirc & \bigcirc
     &\cdot& \bigcirc & \bigcirc & \bigcirc & \bigcirc \\ 
\bigcirc &\cdot& \bigcirc & \bigcirc & \bigcirc & \bigcirc & \bigcirc & \bigcirc
     & \bigcirc &\cdot&\cdot&\cdot&\cdot\\ 
\cdot& \bigcirc &\cdot&\cdot&\cdot&\cdot&\cdot&\cdot&\cdot&\cdot&\cdot&\cdot&\cdot\\ 
\cdot&\cdot&\cdot&\cdot&\cdot&\cdot&\cdot&\cdot&\cdot&\cdot&\cdot&\cdot& \cdot
\end{array}
\]

\[
\gamma^{\natural }= 
\begin{array}{lllllllllllll}
\bigcirc & \bigcirc & \bigcirc & \bigcirc & \bigcirc & \bigcirc & \bigcirc & 
\bigcirc & \bigcirc & \bigcirc & \bigcirc & \bigcirc & \bigcirc \\ 
\cdot& \bigcirc & \bigcirc & \bigcirc & \bigcirc & \bigcirc & \bigcirc & \bigcirc
     & \bigcirc & \bigcirc & \bigcirc & \bigcirc & \bigcirc \\ 
\bigcirc & \bigcirc & \bigcirc & \bigcirc & \bigcirc & \bigcirc & \bigcirc & 
\bigcirc & \bigcirc & \bigcirc & \bigcirc & \bigcirc & \bigcirc \\ 
\cdot& \bigcirc & \bigcirc & \bigcirc & \bigcirc & \bigcirc & \bigcirc & \bigcirc
     &\cdot&\cdot&\cdot&\cdot&\cdot\\ 
\bigcirc &\cdot&\cdot&\cdot&\cdot&\cdot&\cdot&\cdot&\cdot&\cdot&\cdot&\cdot&\cdot\\ 
\cdot&\cdot&\cdot&\cdot&\cdot&\cdot&\cdot&\cdot&\cdot&\cdot&\cdot&\cdot& \cdot
\end{array}
\]

It is possible to prove that $[S(\gamma^{\natural }):D(\beta^{(j)})]=1$ by
applying Rules 2 and 3, but it is tricky to apply Schaper's Theorem without
making a mistake. We therefore prove that $[S(\gamma^{\natural
}):D(\beta^{(j)})]=1$ as follows (recall that $i=2$).

Let $\phi$ be the abacus
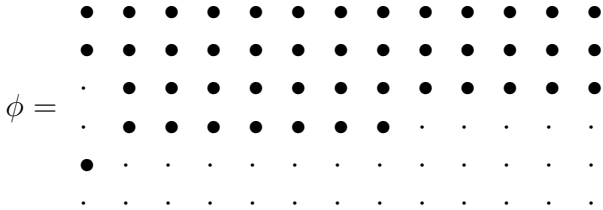
\[
\phi = 
\begin{array}{lllllllllllll}
\bigcirc & \bigcirc & \bigcirc & \bigcirc & \bigcirc & \bigcirc & \bigcirc & 
\bigcirc & \bigcirc & \bigcirc & \bigcirc & \bigcirc & \bigcirc \\ 
\bigcirc & \bigcirc & \bigcirc & \bigcirc & \bigcirc & \bigcirc & \bigcirc & 
\bigcirc & \bigcirc & \bigcirc & \bigcirc & \bigcirc & \bigcirc \\ 
\cdot& \bigcirc & \bigcirc & \bigcirc & \bigcirc & \bigcirc & \bigcirc & \bigcirc
     & \bigcirc & \bigcirc & \bigcirc & \bigcirc & \bigcirc \\ 
\cdot& \bigcirc & \bigcirc & \bigcirc & \bigcirc & \bigcirc & \bigcirc & \bigcirc
     &\cdot&\cdot&\cdot&\cdot&\cdot\\ 
\bigcirc &\cdot&\cdot&\cdot&\cdot&\cdot&\cdot&\cdot&\cdot&\cdot&\cdot&\cdot&\cdot\\ 
\cdot&\cdot&\cdot&\cdot&\cdot&\cdot&\cdot&\cdot&\cdot&\cdot&\cdot&\cdot& \cdot
\end{array}
\]
Using the Littlewood-Richardson rule, we now add a skew $p$-hook to $\phi$ 
in all possible ways to see that:
$$
-\langle 1^{2}\rangle +\langle 1,j\rangle
  -\langle 1,j+1\rangle +\dots \pm \langle 1,p\rangle
  +(-1)^{j}\big(\langle 1,2\rangle -\langle 1,3\rangle
  +\dots\pm \langle 1,j-1\rangle\big)+\langle 1\rangle=0.$$
Now, $\gamma^{\natural }=\langle
1^{2}\rangle $. Therefore, $[S(\gamma^{\natural }):D(\beta^{(j)})]$
is equal to the multiplicity of $D(\beta^{(j)})$ in
\[
\langle 1,j\rangle -\langle 1,j+1\rangle 
  +\dots\pm\langle 1,p\rangle 
  +(-1)^{j}\big(\langle 1,2\rangle -\langle 1,3\rangle 
  +\dots\pm\langle 1,j-1\rangle\big)+\langle 1\rangle.
\]
This is equal to the multiplicity of $D(\beta^{(j)})$ in
$\langle 1,j\rangle -\langle1,j+1\rangle+\dots\pm\langle 1,p\rangle$
because $\beta^{(j)}$ does not dominate the other terms (consider the
first two parts). In turn, this multiplicity is equal to 
$[S(\langle 1,j\rangle):D(\beta^{(j)})]$ since 
$\beta^{(j)}$ does not dominate the other terms ($\beta^{(j)}$ and
all the other terms have the same first column, $\beta^{(j)}$ ends in
$j-2$ ones, while $\langle 1,k\rangle $ ends in $k-2$ ones,
for $j\leq k\leq p$). Finally, $[S(\langle
1,j\rangle):D(\beta^{(j)})]=1$ by two applications of Rule 5
followed by the defect $1$ result.

(3) Suppose that $\mu \in \{\gamma^{\ast },\gamma^{(v)}\}$. Note
that $\mu$ has exactly one normal $(i-1)$-node (except if
$i=2,\;j=p+1$ and $\mu =\gamma^{(v)}$), but every 
exceptional partition $\lambda$ has at most $1$ removable
$(i-1)$-node, so we can apply Proposition 2.1 again. Suppose that
$i=2,\;j=p+1$ and $\mu =\gamma^{(v)}$ (here, $i<v<p$). We need only
consider those $\lambda$ where the first part of $\mu$ exceeds the
first part of $\lambda$, and $\lambda$ has a removable $v$-node. 
It is easily to check that there are no such partitions, so we have finished.
\endproof

We remark that a more detailed analysis shows that the part of the
decomposition matrix with the rows and columns indexed by the
exceptional partitions has the following block diagonal form
$$\begin{array}{ccccc}
         * & 0 & \cdots &\cdots& 0\\
         * & * & 0 &\cdots& 0\\
	 0 & * & * &\cdots& 0\\
	 \vdots&\ddots&\ddots&\ddots&\vdots\\
	 0&\cdots&\ddots&*&*\\
	 0&\cdots&\cdots&0&*
\end{array}$$
where the blocks are certain $3\times3$ matrices (with singular
columns omitted) which are labelled by triples
$\{\alpha^?,\beta^?,\gamma^?\}$. The ordering of $3\times3$ blocks is
compatible with the ordering of the partitions given before the
statement of Proposition 4.1. See the Appendix for the case $p=5$.

\bigskip
\noindent\textbf{Some decomposition numbers in Case 3}

\noindent Assume that the $p$-core belongs to Case 3. That is,
runners up to runner $i-1$ contain $3$ beads; runner $i$ contains $4$
beads; runner $i+1$ contains $5$ beads; after this, there are some or
no runners with $5$ beads; after this, there are some or no runners
with $4$ beads; any remaining runners contain $3$ beads. Let runner
$j$ be the first runner with $4$ beads; let runner $k$ be the first
runner after runner~$i$ with $3$ beads.

By Rule 1, we can equate the decomposition number $[S(\lambda):D(\mu)]$
with a decomposition number of weight $3$ in a smaller symmetric group for
all  partitions $\lambda$ and $\mu$ in the block, except for 
 when  $\lambda$ is one of the partitions $\alpha$, $\beta$, $\gamma$ or $\delta$ 
where
\[
\alpha =\langle i^{2},i+1\rangle ,\quad
\beta  =\langle i-1,i,i+1\rangle ,\quad
\gamma =\langle i^{2}\rangle \quad\text{and}\quad
\delta =\langle i,i-1\rangle . 
\]
We call $\alpha ,\beta ,\gamma ,\delta $ the \textit{exceptional partitions}
for Case 3.  The abacus configurations for the exceptional
partitions in case 3 are as follows:
\[
\alpha =
\begin{array}{lllllllllllll}
\bigcirc  & \bigcirc  & \bigcirc  & \bigcirc  & \bigcirc  & \bigcirc  & 
\bigcirc  & \bigcirc  & \bigcirc  & \bigcirc  & \bigcirc  & \bigcirc  & 
\bigcirc  \\ 
\bigcirc  & \bigcirc  & \bigcirc  & \bigcirc  & \bigcirc  & \bigcirc  & 
\bigcirc  & \bigcirc  & \bigcirc  & \bigcirc  & \bigcirc  & \bigcirc  & 
\bigcirc  \\ 
\bigcirc  & \bigcirc  &\cdot& \bigcirc  & \bigcirc  & \bigcirc  & \bigcirc  & 
\bigcirc  & \bigcirc  & \bigcirc  & \bigcirc  & \bigcirc  & \bigcirc  \\ 
\cdot&\cdot& \bigcirc  & \bigcirc  & \bigcirc  & \bigcirc  & \bigcirc  & \bigcirc  & 
\bigcirc  & \bigcirc  &\cdot&\cdot&\cdot\\ 
\cdot&\cdot& \bigcirc  &\cdot& \bigcirc  & \bigcirc  & \bigcirc  &\cdot&\cdot&\cdot&\cdot&\cdot&\cdot\\ 
\cdot&\cdot&\cdot& \bigcirc  &\cdot&\cdot&\cdot&\cdot&\cdot&\cdot&\cdot&\cdot& \cdot
\end{array}
\]
so $\alpha =(2p-k-j+3i+2,(p-k+2i+1)^{j-i-2},p-k+2i,i^{k-i},1^{p-i})$;
\[
\beta =
\begin{array}{lllllllllllll}
\bigcirc  & \bigcirc  & \bigcirc  & \bigcirc  & \bigcirc  & \bigcirc  & 
\bigcirc  & \bigcirc  & \bigcirc  & \bigcirc  & \bigcirc  & \bigcirc  & 
\bigcirc  \\ 
\bigcirc  & \bigcirc  & \bigcirc  & \bigcirc  & \bigcirc  & \bigcirc  & 
\bigcirc  & \bigcirc  & \bigcirc  & \bigcirc  & \bigcirc  & \bigcirc  & 
\bigcirc  \\ 
\bigcirc  &\cdot& \bigcirc  & \bigcirc  & \bigcirc  & \bigcirc  & \bigcirc  & 
\bigcirc  & \bigcirc  & \bigcirc  & \bigcirc  & \bigcirc  & \bigcirc  \\ 
\cdot& \bigcirc  &\cdot& \bigcirc  & \bigcirc  & \bigcirc  & \bigcirc  & \bigcirc  & 
\bigcirc  & \bigcirc  &\cdot&\cdot&\cdot\\ 
\cdot&\cdot& \bigcirc  &\cdot& \bigcirc  & \bigcirc  & \bigcirc  &\cdot&\cdot&\cdot&\cdot&\cdot&\cdot\\ 
\cdot&\cdot&\cdot& \bigcirc  &\cdot&\cdot&\cdot&\cdot&\cdot&\cdot&\cdot&\cdot& \cdot
\end{array}
\]
so $\beta =(2p-k-j+3i+2,(p-k+2i+1)^{j-i-2},p-k+2i,i^{k-i-1},i-1,1^{p-i+1})$;
\[
\gamma = 
\begin{array}{lllllllllllll}
\bigcirc & \bigcirc & \bigcirc & \bigcirc & \bigcirc & \bigcirc & \bigcirc & 
\bigcirc & \bigcirc & \bigcirc & \bigcirc & \bigcirc & \bigcirc \\ 
\bigcirc & \bigcirc & \bigcirc & \bigcirc & \bigcirc & \bigcirc & \bigcirc & 
\bigcirc & \bigcirc & \bigcirc & \bigcirc & \bigcirc & \bigcirc \\ 
\bigcirc & \bigcirc &\cdot& \bigcirc & \bigcirc & \bigcirc & \bigcirc & \bigcirc
     & \bigcirc & \bigcirc & \bigcirc & \bigcirc & \bigcirc \\ 
\cdot&\cdot& \bigcirc & \bigcirc & \bigcirc & \bigcirc & \bigcirc & \bigcirc & 
\bigcirc & \bigcirc &\cdot&\cdot&\cdot\\ 
\cdot&\cdot&\cdot& \bigcirc & \bigcirc & \bigcirc & \bigcirc &\cdot&\cdot&\cdot&\cdot&\cdot&\cdot\\ 
\cdot&\cdot& \bigcirc &\cdot&\cdot&\cdot&\cdot&\cdot&\cdot&\cdot&\cdot&\cdot& \cdot
\end{array}
\]
so $\gamma =(2p-k-j+3i+1,(p-k+2i+1)^{j-i-1},i^{k-i},1^{p-i})$;
\[
\delta = 
\begin{array}{lllllllllllll}
\bigcirc & \bigcirc & \bigcirc & \bigcirc & \bigcirc & \bigcirc & \bigcirc & 
\bigcirc & \bigcirc & \bigcirc & \bigcirc & \bigcirc & \bigcirc \\ 
\bigcirc & \bigcirc & \bigcirc & \bigcirc & \bigcirc & \bigcirc & \bigcirc & 
\bigcirc & \bigcirc & \bigcirc & \bigcirc & \bigcirc & \bigcirc \\ 
\bigcirc &\cdot& \bigcirc & \bigcirc & \bigcirc & \bigcirc & \bigcirc & \bigcirc
     & \bigcirc & \bigcirc & \bigcirc & \bigcirc & \bigcirc \\ 
\cdot& \bigcirc &\cdot& \bigcirc & \bigcirc & \bigcirc & \bigcirc & \bigcirc & 
\bigcirc & \bigcirc &\cdot&\cdot&\cdot\\ 
\cdot&\cdot&\cdot& \bigcirc & \bigcirc & \bigcirc & \bigcirc &\cdot&\cdot&\cdot&\cdot&\cdot&\cdot\\ 
\cdot&\cdot& \bigcirc &\cdot&\cdot&\cdot&\cdot&\cdot&\cdot&\cdot&\cdot&\cdot& \cdot
\end{array}
\]
so $\delta =(2p-k-j+3i+1,(p-k+2i+1)^{j-i-1},i^{k-i-1},i-1,1^{p-i+1})$.

\bigskip

\noindent 4.2 PROPOSITION \textit{The part of the decomposition matrix whose rows and
columns are labelled by $\alpha ,\beta ,\gamma ,\delta $ is:
\[
\begin{array}{l|llll}
      &\alpha&\beta&\gamma&\delta\\\hline
\alpha&1 & \cdot & \cdot & \cdot \\ 
\beta &1 & 1 & \cdot & \cdot \\ 
\gamma&1 & \cdot & 1 & \cdot \\ 
\delta&1 & 1 & 1 & 1
\end{array}
\]
(Omitted entries are zero.)}

\bigskip

\noindent\textit{Proof.} Note that $\alpha \vartriangleright \beta \vartriangleright
\delta $ and $\alpha \vartriangleright \gamma \vartriangleright \delta $ but 
$\beta \ntriangleright \gamma $ and $\gamma \ntriangleright \beta $.
Let
\begin{eqnarray*}
\alpha_{(1)} &=&\beta_{(1)}=(2p-k-j+3i+2,(p-k+2i+1)^{j-i-2},p-k+2i) \\
\gamma_{(1)} &=&\delta_{(1)}=(2p-k-j+3i+1,(p-k+2i+1)^{j-i-1})
\end{eqnarray*}
and  $\alpha_{(2)}=\gamma_{(2)}=(i^{k-i},1^{p-i})$ and
$\beta_{(2)}=\delta_{(2)}=(i^{k-i-1},i-1,1^{p-i+1})$.
Note that
\[
[S(\gamma_{(1)}):D(\alpha_{(1)})]=[S(p-j+i+1,1^{j-i-1}):D(p-j+i+2,1^{j-i-2})]=1, 
\]
by Rule 5 (column removal), and that
\[
[S(\beta_{(2)}):D(\alpha_{(2)})]=[S(i-1,1^{p-i+1}):D(i,1^{p-i})]=1, 
\]
by Rule 4 (row removal). Next,
\[
[S(\lambda):D(\mu)]=[S(\lambda_{(1)}):D(\mu_{(1)})]
                          [S(\lambda_{(2)}):D(\mu_{(2)})] 
\]
for all $\lambda ,\mu \in \{\alpha ,\beta ,\gamma ,\delta \}$ by Rule 4,
again. From this the Proposition follows.
\endproof

\bigskip
\noindent\textbf{Some decomposition numbers in Case 4}

\noindent Assume that the $p$-core belongs to Case 4. That is, runners
up to runner $i-1$ contain $3$ beads; runner $i$ contains $5$ beads;
after this, there are some or no runners with $5$ beads; after this,
there are some or no runners with $4$ beads; any remaining runners
contain $3$ beads. Let runner $j$ be the first runner with $4$ beads;
let runner $k$ be the first runner after runner $i$ with $3$ beads.

By Rule 1, we can equate the decomposition number $[S(\lambda):D(\mu)]$
with a decomposition number of weight $3$ in a smaller symmetric group for
all $\lambda$ and $\mu$ in the block except for 
 when $\lambda$ is one of the partitions $\alpha$, $\beta$, $\gamma$ or $\delta$
where
\[
\alpha =\langle i^{3}\rangle ,\quad
\beta =\langle i^{2},i-1\rangle ,\,\gamma =\langle i-1,i\rangle ,\quad
\delta =\langle i-1\rangle . 
\]
We call $\alpha ,\beta ,\gamma ,\delta $ the \textit{exceptional partitions}
for Case 4.  The abacus configurations for the exceptional partitions in case 4
are as follows:
\[
\alpha = 
\begin{array}{lllllllllllll}
\bigcirc & \bigcirc & \bigcirc & \bigcirc & \bigcirc & \bigcirc & \bigcirc & 
\bigcirc & \bigcirc & \bigcirc & \bigcirc & \bigcirc & \bigcirc \\ 
\bigcirc & \bigcirc & \bigcirc & \bigcirc & \bigcirc & \bigcirc & \bigcirc & 
\bigcirc & \bigcirc & \bigcirc & \bigcirc & \bigcirc & \bigcirc \\ 
\bigcirc & \bigcirc &\cdot& \bigcirc & \bigcirc & \bigcirc & \bigcirc & \bigcirc
     & \bigcirc & \bigcirc & \bigcirc & \bigcirc & \bigcirc \\ 
\cdot&\cdot& \bigcirc & \bigcirc & \bigcirc & \bigcirc & \bigcirc & \bigcirc & 
\bigcirc & \bigcirc &\cdot&\cdot&\cdot\\ 
\cdot&\cdot& \bigcirc & \bigcirc & \bigcirc & \bigcirc & \bigcirc &\cdot&\cdot&\cdot&\cdot&\cdot& 
\cdot\\ 
\cdot&\cdot& \bigcirc &\cdot&\cdot&\cdot&\cdot&\cdot&\cdot&\cdot&\cdot&\cdot& \cdot
\end{array}
\]
so $\alpha =(2p-k-j+3i,(p-k+2i)^{j-i},i^{k-i},1^{p-i})$;
\[
\beta = 
\begin{array}{lllllllllllll}
\bigcirc & \bigcirc & \bigcirc & \bigcirc & \bigcirc & \bigcirc & \bigcirc & 
\bigcirc & \bigcirc & \bigcirc & \bigcirc & \bigcirc & \bigcirc \\ 
\bigcirc & \bigcirc & \bigcirc & \bigcirc & \bigcirc & \bigcirc & \bigcirc & 
\bigcirc & \bigcirc & \bigcirc & \bigcirc & \bigcirc & \bigcirc \\ 
\bigcirc &\cdot& \bigcirc & \bigcirc & \bigcirc & \bigcirc & \bigcirc & \bigcirc
     & \bigcirc & \bigcirc & \bigcirc & \bigcirc & \bigcirc \\ 
\cdot& \bigcirc &\cdot& \bigcirc & \bigcirc & \bigcirc & \bigcirc & \bigcirc & 
\bigcirc & \bigcirc &\cdot&\cdot&\cdot\\ 
\cdot&\cdot& \bigcirc & \bigcirc & \bigcirc & \bigcirc & \bigcirc &\cdot&\cdot&\cdot&\cdot&\cdot& 
\cdot\\ 
\cdot&\cdot& \bigcirc &\cdot&\cdot&\cdot&\cdot&\cdot&\cdot&\cdot&\cdot&\cdot& \cdot
\end{array}
\]
so $\beta =(2p-k-j+3i,(p-k+2i)^{j-i},i^{k-i-1},i-1,1^{p-i+1})$;
\[
\gamma = 
\begin{array}{lllllllllllll}
\bigcirc & \bigcirc & \bigcirc & \bigcirc & \bigcirc & \bigcirc & \bigcirc & 
\bigcirc & \bigcirc & \bigcirc & \bigcirc & \bigcirc & \bigcirc \\ 
\bigcirc & \bigcirc & \bigcirc & \bigcirc & \bigcirc & \bigcirc & \bigcirc & 
\bigcirc & \bigcirc & \bigcirc & \bigcirc & \bigcirc & \bigcirc \\ 
\bigcirc &\cdot& \bigcirc & \bigcirc & \bigcirc & \bigcirc & \bigcirc & \bigcirc
     & \bigcirc & \bigcirc & \bigcirc & \bigcirc & \bigcirc \\ 
\cdot&\cdot& \bigcirc & \bigcirc & \bigcirc & \bigcirc & \bigcirc & \bigcirc & 
\bigcirc & \bigcirc &\cdot&\cdot&\cdot\\ 
\cdot& \bigcirc &\cdot& \bigcirc & \bigcirc & \bigcirc & \bigcirc &\cdot&\cdot&\cdot&\cdot&\cdot& 
\cdot\\ 
\cdot&\cdot& \bigcirc &\cdot&\cdot&\cdot&\cdot&\cdot&\cdot&\cdot&\cdot&\cdot& \cdot
\end{array}
\]
so $\gamma =(2p-k-j+3i,(p-k+2i)^{j-i-1},p-k+2i-1,i^{k-i},1^{p-i+1})$;
\[
\delta = 
\begin{array}{lllllllllllll}
\bigcirc & \bigcirc & \bigcirc & \bigcirc & \bigcirc & \bigcirc & \bigcirc & 
\bigcirc & \bigcirc & \bigcirc & \bigcirc & \bigcirc & \bigcirc \\ 
\bigcirc & \bigcirc & \bigcirc & \bigcirc & \bigcirc & \bigcirc & \bigcirc & 
\bigcirc & \bigcirc & \bigcirc & \bigcirc & \bigcirc & \bigcirc \\ 
\bigcirc &\cdot& \bigcirc & \bigcirc & \bigcirc & \bigcirc & \bigcirc & \bigcirc
     & \bigcirc & \bigcirc & \bigcirc & \bigcirc & \bigcirc \\ 
\cdot&\cdot& \bigcirc & \bigcirc & \bigcirc & \bigcirc & \bigcirc & \bigcirc & 
\bigcirc & \bigcirc &\cdot&\cdot&\cdot\\ 
\cdot&\cdot& \bigcirc & \bigcirc & \bigcirc & \bigcirc & \bigcirc &\cdot&\cdot&\cdot&\cdot&\cdot& 
\cdot\\ 
\cdot& \bigcirc &\cdot&\cdot&\cdot&\cdot&\cdot&\cdot&\cdot&\cdot&\cdot&\cdot& \cdot
\end{array}
\]
so $\delta =(2p-k-j+3i-1,(p-k+2i)^{j-i},i^{k-i},1^{p-i+1})$.

\bigskip

\noindent 4.3 PROPOSITION \textit{The part of the decomposition matrix whose rows and
columns are labelled by $\alpha ,\beta ,\gamma ,\delta $ is:
\[
\begin{array}{l|llll}
      &\alpha&\beta&\gamma&\delta\\\hline
\alpha&1 & \cdot & \cdot & \cdot \\ 
\beta &1 & 1 & \cdot & \cdot \\ 
\gamma&1 & 1 & 1 & \cdot \\ 
\delta&1 & 1 & 1 & 1
\end{array}
\]
(Omitted entries are zero.)}

\noindent\textit{Proof.} Note that $\alpha \vartriangleright \beta \vartriangleright
\gamma \vartriangleright \delta $. First, each of $\alpha ,\beta ,\gamma $
and $\delta $ has exactly $3$ removable $i$-nodes. Also, the $3$ removable 
$ i $-nodes in $\alpha $ are normal. Therefore, for 
$\lambda \in \{\alpha ,\beta ,\gamma ,\delta \}$ we have 
$[S(\lambda):D(\alpha)]=[S(\bar{\lambda}):D(\bar{\alpha})]$ as in
Proposition 2.1.  However, here 
$\bar{\alpha}=\bar{\beta }=\bar{\gamma}=\bar{\delta}$, so 
$[S(\lambda):D(\alpha)]=1$.

Next, remove the first column from $\beta ,\gamma ,\delta $; see Note (a)
which follows Rule 5. We can now apply a similar argument to the above,
using $2$ removable $i$-nodes, to deduce that $[S(\lambda):D(\beta)]=1$
for $\lambda \in \{\beta ,\gamma ,\delta \}$.

Finally, remove the first $i$ columns from $\gamma ,\delta $, 
using Rule 5, and use the one remaining removable $i$-node to deduce that
$[S(\lambda):D(\gamma)]=1$, for $\lambda \in \{\gamma ,\delta \}$. The
proof of the Proposition is now complete.
\endproof

Of course, Propositions 4.1, 4.2 and 4.3 hardly scratch the surface of the
problem of calculating the decomposition numbers, since they evaluate 
$[S(\lambda):D(\mu)]$ only when $\lambda$ and $\mu$ are both exceptional.
It is still necessary to calculate $[S(\lambda):D(\mu)]$ when $\lambda$ is
exceptional and $\mu$ is arbitrary. In Case 4, for example, these answers
depend upon the values of $i,j$ and $k$, and 
there are a very large number of separate cases that have to be
considered.

\bigskip

\[
\text{5. THE CASE\ }p=5 
\]

\bigskip

We have written a computer program,  using the \textsf{GAP} package
\textsc{Specht} [\textbf{M1}], to apply Rules 1-7 to find all the
decomposition numbers when $p=5$ and $w=3$. All the decomposition
numbers were determined once we had assumed L\"{u}beck and
M\"{u}ller's result that $[S(8^{2},4,1):D(12,9)]=1$.  Consequently,
 at least in principle, all of the decomposition numbers of the
symmetric groups for blocks of weight $3$ in characteristic $5$ are
now known. The decomposition matrices for Cases 2--4 when $p=5$ are
given in the appendix.

In outline, the program first finds all the partitions $\nu $ which
dominate the last of the exceptional partitions and then uses
Rules~1 and 2 to find the decomposition numbers $[S(\nu):D(\mu
)]$ whenever $\nu $ is not exceptional. For the exceptional partitions
$\lambda$, the program applies Schaper's Theorem (Rule 2); this often
determines the decomposition numbers $[S(\lambda):D(\mu)]$. If this
decomposition number is not determined  then Rule 2 gives us an
integer $m>1$ such that $m\geq [S(\lambda):D(\mu)]\geq 1$. The
program next checks to see whether the answer is given by one of Rules
3--6. Finally, as a last resort, the program tries to apply Rule 7 in
order to show that $[S(\lambda):D(\mu)]=1$.  The program also does
parallel computations with two different primes which it uses, along
with Rule~1, to check the consistency of its calculations (compare
Conjecture 3.1). 

Finally, in order to check our calculations we compared the matrices 
that we computed with the decomposition matrices of the corresponding
Hecke algebra of type $A$ [\textbf{M2}] at a complex $p^{\text{th}}$
root of unity --- which are known by the LLT algorithm [\textbf{LLT}].  
Since we were able to compute these decomposition
numbers using only Rules 1--7 (and L\"{u}beck and J. M\"{u}ller's
result) these two sets of decomposition matrices should agree because
$p>w$ (this affects only Rule 2). In all cases the symmetric group and
Hecke algebra decomposition multiplicities were the same.

Here is a small example of the technique in action.

\bigskip

\noindent\textbf{Example} Suppose that we are in Case 4, with
$i=p,\,j=k=p+1$. Thus the $p$-core has the following abacus:
\[
\begin{array}{ccccc}
\bigcirc & \bigcirc & \bigcirc & \bigcirc & \bigcirc \\ 
\bigcirc & \bigcirc & \bigcirc & \bigcirc & \bigcirc \\ 
\bigcirc & \bigcirc & \bigcirc & \bigcirc & \bigcirc \\ 
\cdot&\cdot&\cdot&\cdot& \bigcirc \\ 
\cdot&\cdot&\cdot&\cdot& \bigcirc
\end{array}
\]
Then the exceptional partitions for this core are: 
$\alpha =\langle p^{3}\rangle$, 
$\beta =\langle p^{2},p-1\rangle$,
$\gamma =\langle p-1,p\rangle$
and $\delta=\langle p-1\rangle$. We will show that the non-zero
decomposition numbers for these exceptional partitions are as
follows:
\[
\begin{array}{ll|ccccccc}
&&\langle p\rangle 
& \langle p^{2}\rangle 
&  \langle p,p-1\rangle &\alpha&\beta&\gamma&\delta\\\hline
         & \langle p\rangle         & 1 &    &    &    &    &    &  \\ 
         & \langle p^{2}\rangle     &\cdot& 1   &    &    &    &    &  \\ 
         & \langle p,p-1\rangle     & 1   & 1   & 1 &    &    &    &  \\ 
\alpha = & \langle p^{3}\rangle     &\cdot&\cdot&\cdot& 1 &    &    &  \\ 
\beta =  & \langle p^{2},p-1\rangle &\cdot& 1   &\cdot& 1 & 1 &    &  \\ 
\gamma = & \langle p-1,p\rangle     &\cdot& 1   & 1   & 1 & 1 & 1 &  \\ 
\delta = & \langle p-1\rangle       &\cdot&\cdot&\cdot& 1 & 1 & 1 & 1
\end{array}
\]

First, one readily checks that the partitions  which index the rows
of this matrix are precisely the partitions $\nu$ such that $\nu$ has
the same $p$-core as $\delta$ and $\nu\trianglerighteq\delta$.

Suppose that $\nu ,\mu \in \{\langle p\rangle ,\langle
p^{2}\rangle ,\langle p,p-1\rangle \}$. Then $p-1$
applications of Rule 1 allow us to equate $[S(\nu):D(\mu)]$ with a
decomposition number in Case 2 with $i=p,\,j=k=p+1$. In practice, though, it
is much easier to apply Rule 2 (Schaper's Theorem) to evaluate $[S(\nu
):D(\mu)]$.

Now consider the exceptional partitions $\alpha ,\beta ,\gamma ,\delta $.

Rule 2 immediately implies that $S(\alpha)$ is irreducible.

 Next, applying Schaper's Theorem to $\beta$ gives the following linear
combination of Specht modules:
$S(\langle p^{2}\rangle)+S(\alpha)$. From what we
have already deduced, this is equal to 
$D(\langle p^{2}\rangle)+D(\alpha)$. Rule 2 now gives us
the row of the matrix which is labelled by~$\beta $.

Similarly, applying Schaper's Theorem to $\gamma$ gives
$$-S(\langle p\rangle)+S(\langle p,p-1\rangle)+S(\alpha)
    +S(\beta)
    =
2D(\langle p^{2}\rangle)+D(\langle p,p-1\rangle
)+2D(\alpha)+D(\beta). 
$$
 By Rule 2, $[S(\gamma):D(\langle p,p-1\rangle)]=[S(\gamma):D(\beta)]=1$
and we also know that
$$ 2 \geq [S(\gamma):D(\langle p^{2}\rangle)]\geq 1
       \quad \text{and}\quad  2 \geq [S(\gamma):D(\alpha)]\geq 1.$$
Now, $\gamma $ and $\langle p^{2}\rangle $ have the same first
part so 
 $[S(\gamma):D(\langle p^{2}\rangle)]
      =[S(\gamma_{(2)}):D(\langle p^{2}\rangle_{(2)})]$ 
by Rule 4. But $\gamma_{(2)}$ and $\langle p^{2}\rangle_{(2)}$
belong to a block of weight 2, so 
$[S(\gamma_{(2)}):D(\langle p^{2}\rangle_{(2)})]\leq 1$. 
Hence $[S(\gamma):D(\langle p^{2}\rangle)]=1$.

If $[S(\gamma):D(\alpha)]=2,$ then

\[
D(\gamma)=S(\langle p-1,p\rangle)-S(\langle
p^{2},p-1\rangle)-S(\langle p,p-1\rangle)+S(\langle
p\rangle). 
\]
If we $p^{2}(p-1)^{2}\dots2^{2}1^{2}$ restrict this, as in Rule 6, we do
not obtain a module (an irreducible module occurs with negative
multiplicity). This contradiction implies that $[S(\gamma):D(\alpha)]=1,$
and all the decomposition numbers for $S(\gamma)$ are now known.

Finally, we apply Schaper's Theorem to $\delta $. This gives
\[
S(\langle p\rangle)-S(\langle p^{2}\rangle)
    -S(\langle p,p-1\rangle)+S(\alpha)+S(\beta)+S(\gamma) 
       =3D(\alpha)+2D(\beta)+D(\gamma).
\]
Thus, the only decomposition numbers for $S(\delta)$ which are still in doubt are 
$[S(\delta):D(\alpha)]$ and $[S(\delta):D(\beta)]$.

We apply Rule 3, with the Kleshchev sequence $p^{3}$ (which leads to a block
of weight $0$) to conclude that $[S(\delta):D(\alpha)]\leq 1$. Hence, 
$[S(\delta):D(\alpha)]=1$.

Now, $\delta $ and $\beta $ have the same first column, so by Rule 5,
\[
[S(\delta):D(\beta)]=[S(\delta^{(1)}:D(\beta^{(1)})]. 
\]
But $\delta^{(1)}$ and $\beta^{(1)}$ belong to a block of weight 2, so 
$[S(\delta^{(1)}:D(\beta^{(1)})]\leq 1$. Hence $[S(\delta):D(\gamma)]=1$.

We have now completed the example.

\bigskip

\noindent\textbf{Notes}

(a) Some of the decomposition numbers in the last example were computed in
different ways in Proposition 4.3. The method in the example uses only Rules
1-7 (in fact we used all the Rules except Rule 7).

(b) The arguments used in the example apply equally well for any $p>3$. As a
consequence of many other instances of this phenomenon, we were led to
formulate Conjecture 3.1.

\bigskip

\[
\text{ACKNOWLEDGMENT} 
\]

This research was supported, in part, by EPSRC Grant GR/S06639/01.

\bigskip 
\[
\text{Appendix: DECOMPOSITION NUMBERS OF WEIGHT 3 IN CHARACTERISTIC 5} 
\]

\noindent  In this appendix we list the non-zero entries in the rows
indexed by exceptional partitions for all of the decomposition
matrices in Cases 2--4 when $p=5$. These matrices, combined with the
results of this paper (specifically Rule 1), determine the
decomposition matrices for all blocks of weight $3$ for all symmetric
groups when $p=5$. 

We remark that we have also used our program to calculate the
decomposition numbers in Cases 1--4 when $p=7$. This calculation took
over one month to complete, on a reasonably fast computer. We were
unable to determine whether the following two decomposition numbers
are equal $1$ or $2$:
$$\begin{array}{l@{\qquad}l}
\text{Case 2:  $(i,j)=(p-1, p+1)$} &  \text{Case 2:  $(i,j)=(p-2, p+1)$}\\\relax
[S(\langle (p-1)^2,p\rangle, D(\langle p-1,p\rangle)] &
  [S(\langle (p-1)^2,p-1\rangle, D(\langle p,p-1,p-2\rangle)]\\
  \text{Core: $(5^2)=\big((p-2)^2\big)$} & \text{Core: $(4^3)=\big((p-2)^3\big)$}
\end{array}$$
Conjecture 3.1 and our calculations for $p=5$ imply that 
of these decomposition numbers should both be equal to $1$. \textit{Assuming
this}, we were able to compute all of the remaining decomposition
numbers for Cases 2--4 when $p=7$. We again found that
$[S(\lambda):D(\mu)]\le1$ in all cases.

\def\thesubsection{A\arabic{subsection}}
\def\L#1{\begin{sideways}$#1$\end{sideways}}

\small
\setlength\arraycolsep{3pt}
\subsection{Case 2: $(i,j)=(p, p+1)$}
\begin{sideways}
$\begin{array}{l|*{23}c}
  &\L{\langle p\rangle}
  &\L{\langle p-1,p\rangle}
  &\L{\langle p^{2},p-1\rangle}
  &\L{\langle p,p-2\rangle}
  &\L{\langle p-2\rangle}
  &\L{\langle p-2,p\rangle}
  &\L{\langle p,p-3\rangle}
  &\L{\langle p-3\rangle}
  &\L{\langle p-3,p\rangle}
  &\L{\langle p-4\rangle}
  &\L{\langle p-4,p\rangle}
  &\L{\langle p^{2}\rangle}
  &\L{\langle p,p-1\rangle}
  &\L{\langle p-1\rangle}
  &\L{\langle p^{2},p-2\rangle}
  &\L{\langle p,p-1,p-2\rangle}
  &\L{\langle p-1,p-2\rangle}
  &\L{\langle p^{2},p-3\rangle}
  &\L{\langle p,p-1,p-3\rangle}
  &\L{\langle p-1,p-3\rangle}
  &\L{\langle p^{2},p-4\rangle}
  &\L{\langle p,p-1,p-4\rangle}
  &\L{\langle p-1,p-4\rangle}\\\toprule
\langle p^{2}\rangle&.&.&.&.&.&.&.&.&.&.&.&1&.&.&.&.&.&.&.&.&.&.&.\\
\langle p,p-1\rangle&1&.&.&.&.&.&.&.&.&.&.&1&1&.&.&.&.&.&.&.&.&.&.\\
\langle p-1\rangle&.&.&.&.&.&.&.&.&.&.&.&1&1&1&.&.&.&.&.&.&.&.&.\\
\langle p^{2},p-2\rangle&1&1&1&1&1&1&.&.&.&.&.&1&1&1&1&.&.&.&.&.&.&.&.\\
\langle p,p-1,p-2\rangle&1&.&.&.&.&.&.&.&.&.&.&1&.&.&1&1&.&.&.&.&.&.&.\\
\langle p-1,p-2\rangle&1&.&1&1&.&1&.&.&.&.&.&1&1&.&1&1&1&.&.&.&.&.&.\\
\langle p^{2},p-3\rangle&1&.&.&1&1&1&1&1&1&.&.&.&.&.&1&.&.&1&.&.&.&.&.\\
\langle p,p-1,p-3\rangle&.&.&1&.&.&.&.&.&.&.&.&.&.&.&1&1&.&1&1&.&.&.&.\\
\langle p-1,p-3\rangle&1&1&1&1&.&1&1&.&1&.&.&.&.&.&1&1&1&1&1&1&.&.&.\\
\langle p^{2},p-4\rangle&.&.&.&.&.&.&1&1&1&1&1&.&.&.&.&.&.&1&.&.&1&.&.\\
\langle p,p-1,p-4\rangle&.&.&.&.&.&.&.&.&.&.&.&.&.&.&.&.&.&1&1&.&1&1&.\\
\langle p-1,p-4\rangle&.&.&.&.&.&.&1&.&1&.&1&.&.&.&.&.&.&1&1&1&1&1&1\\
\langle p^{3}\rangle&.&.&.&.&.&.&.&.&.&1&1&.&.&.&.&.&.&.&.&.&1&.&.\\
\langle (p-1)^{2},p\rangle&.&.&.&.&.&.&.&.&.&.&.&.&.&.&.&.&.&.&.&.&1&1&.\\
\langle (p-1)^{2}\rangle&.&.&.&.&.&.&.&.&.&.&1&.&.&.&.&.&.&.&.&.&1&1&1\\
\end{array}$
\end{sideways}

\subsection{Case 2: $(i,j)=(p-1, p)$}
\begin{sideways}
$\begin{array}{l|*{23}c}
  &\L{\langle p-1\rangle}
  &\L{\langle p\rangle}
  &\L{\langle p-1,p\rangle}
  &\L{\langle p,p-1\rangle}
  &\L{\langle p,p-2\rangle}
  &\L{\langle p-2,p-1\rangle}
  &\L{\langle (p-1)^{2},p-2\rangle}
  &\L{\langle p,p-3\rangle}
  &\L{\langle p,p-1,p-3\rangle}
  &\L{\langle p,p-4\rangle}
  &\L{\langle p,p-1,p-4\rangle}
  &\L{\langle (p-1)^{2}\rangle}
  &\L{\langle p-1,p-2\rangle}
  &\L{\langle p-2\rangle}
  &\L{\langle (p-1)^{2},p\rangle}
  &\L{\langle p,p-1,p-2\rangle}
  &\L{\langle p-2,p\rangle}
  &\L{\langle (p-1)^{2},p-3\rangle}
  &\L{\langle p-1,p-2,p-3\rangle}
  &\L{\langle p-2,p-3\rangle}
  &\L{\langle (p-1)^{2},p-4\rangle}
  &\L{\langle p-1,p-2,p-4\rangle}
  &\L{\langle p-2,p-4\rangle}\\\toprule
\langle (p-1)^{2}\rangle&1&1&1&.&.&.&.&.&.&.&.&1&.&.&.&.&.&.&.&.&.&.&.\\
\langle p-1,p-2\rangle&.&.&.&.&.&.&.&.&.&.&.&1&1&.&.&.&.&.&.&.&.&.&.\\
\langle p-2\rangle&1&.&1&.&.&.&.&.&.&.&.&1&1&1&.&.&.&.&.&.&.&.&.\\
\langle (p-1)^{2},p\rangle&.&.&1&1&.&.&.&.&.&.&.&1&.&.&1&.&.&.&.&.&.&.&.\\
\langle p,p-1,p-2\rangle&1&1&1&1&1&1&.&.&.&.&.&1&1&1&1&1&.&.&.&.&.&.&.\\
\langle p-2,p\rangle&1&.&.&.&.&.&.&.&.&.&.&1&1&.&1&1&1&.&.&.&.&.&.\\
\langle (p-1)^{2},p-3\rangle&.&.&.&1&1&.&1&1&1&.&.&.&.&.&1&1&.&1&.&.&.&.&.\\
\langle p-1,p-2,p-3\rangle&.&.&.&.&.&.&.&.&.&.&.&.&.&.&1&.&.&1&1&.&.&.&.\\
\langle p-2,p-3\rangle&.&.&.&1&.&1&1&.&1&.&.&.&.&.&1&1&1&1&1&1&.&.&.\\
\langle (p-1)^{2},p-4\rangle&.&.&.&.&.&.&.&1&1&1&1&.&.&.&.&.&.&1&.&.&1&.&.\\
\langle p-1,p-2,p-4\rangle&.&.&.&.&.&.&1&.&.&.&.&.&.&.&.&.&.&1&1&.&1&1&.\\
\langle p-2,p-4\rangle&.&.&.&.&.&.&1&.&1&.&1&.&.&.&.&.&.&1&1&1&1&1&1\\
\langle (p-1)^{3}\rangle&.&.&.&.&.&.&.&.&.&1&1&.&.&.&.&.&.&.&.&.&1&.&.\\
\langle (p-2)^{2},p-1\rangle&.&.&.&.&.&.&.&.&.&.&.&.&.&.&.&.&.&.&.&.&1&1&.\\
\langle (p-2)^{2}\rangle&.&.&.&.&.&.&.&.&.&.&1&.&.&.&.&.&.&.&.&.&1&1&1\\
\end{array}$
\end{sideways}

\subsection{Case 2: $(i,j)=(p-2, p-1)$}
\begin{sideways}
$\begin{array}{l|*{26}c}
  &\L{\langle p-2\rangle}
  &\L{\langle p\rangle}
  &\L{\langle p,p-2\rangle}
  &\L{\langle p-1\rangle}
  &\L{\langle p-2,p-1\rangle}
  &\L{\langle p,p-1\rangle}
  &\L{\langle p-1,p-2\rangle}
  &\L{\langle p-1,p\rangle}
  &\L{\langle p,p-1,p-2\rangle}
  &\L{\langle p,p-3\rangle}
  &\L{\langle p-1,p-3\rangle}
  &\L{\langle p-3,p-2\rangle}
  &\L{\langle (p-2)^{2},p-3\rangle}
  &\L{\langle p-1,p-4\rangle}
  &\L{\langle p-1,p-2,p-4\rangle}
  &\L{\langle (p-2)^{2}\rangle}
  &\L{\langle p-2,p-3\rangle}
  &\L{\langle p-3\rangle}
  &\L{\langle (p-2)^{2},p\rangle}
  &\L{\langle p,p-2,p-3\rangle}
  &\L{\langle p-3,p\rangle}
  &\L{\langle (p-2)^{2},p-1\rangle}
  &\L{\langle p-1,p-2,p-3\rangle}
  &\L{\langle p-3,p-1\rangle}
  &\L{\langle (p-2)^{2},p-4\rangle}
  &\L{\langle p-3,p-4\rangle}\\\toprule
\langle (p-2)^{2}\rangle&.&.&.&1&1&.&.&.&.&.&.&.&.&.&.&1&.&.&.&.&.&.&.&.&.&.\\
\langle p-2,p-3\rangle&.&.&.&.&.&.&.&.&.&.&.&.&.&.&.&1&1&.&.&.&.&.&.&.&.&.\\
\langle p-3\rangle&.&.&.&.&1&.&.&.&.&.&.&.&.&.&.&1&1&1&.&.&.&.&.&.&.&.\\
\langle (p-2)^{2},p\rangle&1&1&.&1&1&1&1&1&1&.&.&.&.&.&.&1&.&.&1&.&.&.&.&.&.&.\\
\langle p,p-2,p-3\rangle&1&.&1&.&1&1&1&.&.&1&.&1&.&.&.&1&1&1&1&1&.&.&.&.&.&.\\
\langle p-3,p\rangle&.&.&.&.&.&.&.&.&1&.&.&.&.&.&.&1&1&.&1&1&1&.&.&.&.&.\\
\langle (p-2)^{2},p-1\rangle&1&1&.&.&.&1&.&.&.&.&.&.&.&.&.&.&.&.&1&.&.&1&.&.&.&.\\
\langle p-1,p-2,p-3\rangle&.&.&1&.&.&1&1&1&1&1&1&.&.&.&.&.&.&.&1&1&.&1&1&.&.&.\\
\langle p-3,p-1\rangle&1&.&.&.&.&.&1&.&1&.&.&1&.&.&.&.&.&.&1&1&1&1&1&1&.&.\\
\langle (p-2)^{2},p-4\rangle&.&.&.&.&.&.&.&.&.&.&1&.&1&1&1&.&.&.&.&.&.&1&1&.&1&.\\
\langle p-2,p-3,p-4\rangle&.&.&.&.&.&.&.&.&.&.&.&.&.&.&.&.&.&.&.&.&.&1&.&.&1&.\\
\langle p-3,p-4\rangle&.&.&.&.&.&.&.&.&.&.&.&.&1&.&1&.&.&.&.&.&.&1&1&1&1&1\\
\langle (p-2)^{3}\rangle&.&.&.&.&.&.&.&.&.&.&.&.&.&1&1&.&.&.&.&.&.&.&.&.&1&.\\
\langle (p-3)^{2},p-2\rangle&.&.&.&.&.&.&.&.&.&.&.&.&1&.&.&.&.&.&.&.&.&.&.&.&1&.\\
\langle (p-3)^{2}\rangle&.&.&.&.&.&.&.&.&.&.&.&.&1&.&1&.&.&.&.&.&.&.&.&.&1&1\\
\end{array}$
\end{sideways}

\subsection{Case 2: $(i,j)=(p-3, p-2)$}
\begin{sideways}
$\begin{array}{l|*{28}c}
  &\L{\langle p,p-3\rangle}
  &\L{\langle p-1,p-3\rangle}
  &\L{\langle p-1,p\rangle}
  &\L{\langle p,p-1,p-3\rangle}
  &\L{\langle p-2\rangle}
  &\L{\langle p-3,p-2\rangle}
  &\L{\langle p,p-2\rangle}
  &\L{\langle p-1,p-2\rangle}
  &\L{\langle p-2,p-3\rangle}
  &\L{\langle p-2,p\rangle}
  &\L{\langle p-2,p-1\rangle}
  &\L{\langle p,p-2,p-3\rangle}
  &\L{\langle p-1,p-2,p-3\rangle}
  &\L{\langle p,p-4\rangle}
  &\L{\langle p-1,p-4\rangle}
  &\L{\langle p-2,p-4\rangle}
  &\L{\langle p-4,p-3\rangle}
  &\L{\langle (p-3)^{2}\rangle}
  &\L{\langle p-4\rangle}
  &\L{\langle (p-3)^{2},p\rangle}
  &\L{\langle p,p-3,p-4\rangle}
  &\L{\langle p-4,p\rangle}
  &\L{\langle (p-3)^{2},p-1\rangle}
  &\L{\langle p-1,p-3,p-4\rangle}
  &\L{\langle p-4,p-1\rangle}
  &\L{\langle (p-3)^{2},p-2\rangle}
  &\L{\langle p-2,p-3,p-4\rangle}
  &\L{\langle p-4,p-2\rangle}\\\toprule
\langle (p-3)^{2}\rangle&.&.&.&.&1&1&.&.&.&.&.&.&.&.&.&.&.&1&.&.&.&.&.&.&.&.&.&.\\
\langle p-3,p-4\rangle&.&.&.&.&.&.&.&.&.&.&.&.&.&.&.&.&.&1&.&.&.&.&.&.&.&.&.&.\\
\langle p-4\rangle&.&.&.&.&.&1&.&.&.&.&.&.&.&.&.&.&.&1&1&.&.&.&.&.&.&.&.&.\\
\langle (p-3)^{2},p\rangle&1&.&.&.&1&1&1&.&1&1&.&1&.&.&.&.&.&1&.&1&.&.&.&.&.&.&.&.\\
\langle p,p-3,p-4\rangle&.&.&.&.&.&1&1&.&1&.&.&.&.&1&.&.&1&1&1&1&1&.&.&.&.&.&.&.\\
\langle p-4,p\rangle&1&.&.&.&.&.&.&.&.&.&.&1&.&.&.&.&.&1&.&1&1&1&.&.&.&.&.&.\\
\langle (p-3)^{2},p-1\rangle&1&1&1&1&.&.&1&1&1&1&1&1&1&.&.&.&.&.&.&1&.&.&1&.&.&.&.&.\\
\langle p-1,p-3,p-4\rangle&.&.&.&.&.&.&1&1&.&.&.&.&.&1&1&.&.&.&.&1&1&.&1&1&.&.&.&.\\
\langle p-4,p-1\rangle&1&1&.&1&.&.&.&.&1&.&.&1&1&.&.&.&1&.&.&1&1&1&1&1&1&.&.&.\\
\langle (p-3)^{2},p-2\rangle&.&1&1&1&.&.&.&1&.&.&.&.&.&.&.&.&.&.&.&.&.&.&1&.&.&1&.&.\\
\langle p-2,p-3,p-4\rangle&.&.&.&.&.&.&.&1&.&.&1&.&1&.&1&1&.&.&.&.&.&.&1&1&.&1&1&.\\
\langle p-4,p-2\rangle&.&1&.&1&.&.&.&.&.&.&.&.&1&.&.&.&.&.&.&.&.&.&1&1&1&1&1&1\\
\langle (p-3)^{3}\rangle&.&.&.&.&.&.&.&.&.&.&.&.&.&.&.&1&.&.&.&.&.&.&.&.&.&1&1&.\\
\langle (p-4)^{2},p-3\rangle&.&.&.&.&.&.&.&.&.&.&.&.&.&.&.&.&.&.&.&.&.&.&.&.&.&1&.&.\\
\langle (p-4)^{2}\rangle&.&.&.&.&.&.&.&.&.&.&.&.&.&.&.&.&.&.&.&.&.&.&.&.&.&1&1&1\\
\end{array}$
\end{sideways}

\subsection{Case 2: $(i,j)=(p-1, p+1)$}
\begin{sideways}
$\begin{array}{l|*{31}c}
  &\L{\langle p\rangle}
  &\L{\langle p-1\rangle}
  &\L{\langle p,p-1\rangle}
  &\L{\langle p-1,p\rangle}
  &\L{\langle p^{2}\rangle}
  &\L{\langle p^{2},p-1\rangle}
  &\L{\langle p,p-2\rangle}
  &\L{\langle p^{2},p-2\rangle}
  &\L{\langle p-2,p-1\rangle}
  &\L{\langle (p-1)^{2},p-2\rangle}
  &\L{\langle p,p-1,p-3\rangle}
  &\L{\langle p-3,p\rangle}
  &\L{\langle p^{2},p-3\rangle}
  &\L{\langle p,p-2,p-3\rangle}
  &\L{\langle p-3,p-1\rangle}
  &\L{\langle p-4,p\rangle}
  &\L{\langle p^{2},p-4\rangle}
  &\L{\langle p,p-2,p-4\rangle}
  &\L{\langle p-4,p-1\rangle}
  &\L{\langle (p-1)^{2}\rangle}
  &\L{\langle p-1,p-2\rangle}
  &\L{\langle p-2\rangle}
  &\L{\langle (p-1)^{2},p\rangle}
  &\L{\langle p,p-1,p-2\rangle}
  &\L{\langle p-2,p\rangle}
  &\L{\langle (p-1)^{2},p-3\rangle}
  &\L{\langle p-1,p-2,p-3\rangle}
  &\L{\langle p-2,p-3\rangle}
  &\L{\langle (p-1)^{2},p-4\rangle}
  &\L{\langle p-1,p-2,p-4\rangle}
  &\L{\langle p-2,p-4\rangle}\\\toprule
\langle (p-1)^{2}\rangle&1&1&1&1&1&.&.&.&.&.&.&.&.&.&.&.&.&.&.&1&.&.&.&.&.&.&.&.&.&.&.\\
\langle p-1,p-2\rangle&.&.&.&.&1&.&1&.&.&.&.&.&.&.&.&.&.&.&.&1&1&.&.&.&.&.&.&.&.&.&.\\
\langle p-2\rangle&.&1&.&1&.&.&.&.&.&.&.&.&.&.&.&.&.&.&.&1&1&1&.&.&.&.&.&.&.&.&.\\
\langle (p-1)^{2},p\rangle&.&1&1&1&1&1&.&.&.&.&.&.&.&.&.&.&.&.&.&1&.&.&1&.&.&.&.&.&.&.&.\\
\langle p,p-1,p-2\rangle&.&.&.&1&.&.&.&.&.&.&.&.&.&.&.&.&.&.&.&1&1&.&1&1&.&.&.&.&.&.&.\\
\langle p-2,p\rangle&1&1&.&1&1&1&1&.&.&.&.&.&.&.&.&.&.&.&.&1&1&1&1&1&1&.&.&.&.&.&.\\
\langle (p-1)^{2},p-3\rangle&.&.&.&1&.&1&.&1&1&1&1&1&1&.&1&.&.&.&.&.&.&.&1&1&1&1&.&.&.&.&.\\
\langle p-1,p-2,p-3\rangle&.&.&.&1&.&1&.&.&.&.&.&.&1&1&.&.&.&.&.&.&.&.&1&.&.&1&1&.&.&.&.\\
\langle p-2,p-3\rangle&.&.&.&1&.&.&.&.&.&1&1&.&.&.&1&.&.&.&.&.&.&.&1&1&.&1&1&1&.&.&.\\
\langle (p-1)^{2},p-4\rangle&.&.&.&1&.&.&.&.&.&.&1&1&1&.&1&1&1&.&1&.&.&.&.&.&.&1&.&.&1&.&.\\
\langle p-1,p-2,p-4\rangle&.&.&.&.&.&.&.&1&.&1&.&.&1&1&.&.&1&1&.&.&.&.&.&.&.&1&1&.&1&1&.\\
\langle p-2,p-4\rangle&.&.&.&1&.&.&.&.&1&1&1&.&.&.&1&.&.&.&1&.&.&.&.&.&.&1&1&1&1&1&1\\
\langle (p-1)^{3}\rangle&.&.&.&.&.&.&.&.&.&.&.&.&.&.&.&1&1&.&1&.&.&.&.&.&.&.&.&.&1&.&.\\
\langle (p-2)^{2},p-1\rangle&.&.&.&.&.&.&.&.&.&.&.&.&.&.&.&.&1&1&.&.&.&.&.&.&.&.&.&.&1&1&.\\
\langle (p-2)^{2}\rangle&.&.&.&.&.&.&.&.&.&.&.&.&.&.&.&.&.&.&1&.&.&.&.&.&.&.&.&.&1&1&1\\
\end{array}$
\end{sideways}

\subsection{Case 2: $(i,j)=(p-2, p)$}
\begin{sideways}
$\begin{array}{l|*{28}c}
  &\L{\langle p-1\rangle}
  &\L{\langle p-2,p-1\rangle}
  &\L{\langle p-1,p\rangle}
  &\L{\langle p-2,p\rangle}
  &\L{\langle p,p-1\rangle}
  &\L{\langle p,p-1,p-2\rangle}
  &\L{\langle (p-1)^{2}\rangle}
  &\L{\langle (p-1)^{2},p-2\rangle}
  &\L{\langle (p-1)^{2},p\rangle}
  &\L{\langle p-1,p-3\rangle}
  &\L{\langle p,p-1,p-3\rangle}
  &\L{\langle (p-1)^{2},p-3\rangle}
  &\L{\langle p-3,p-2\rangle}
  &\L{\langle (p-2)^{2},p-3\rangle}
  &\L{\langle p,p-1,p-4\rangle}
  &\L{\langle (p-1)^{2},p-4\rangle}
  &\L{\langle p,p-2,p-4\rangle}
  &\L{\langle (p-2)^{2}\rangle}
  &\L{\langle p-2,p-3\rangle}
  &\L{\langle p-3\rangle}
  &\L{\langle (p-2)^{2},p-1\rangle}
  &\L{\langle p-1,p-2,p-3\rangle}
  &\L{\langle p-3,p-1\rangle}
  &\L{\langle (p-2)^{2},p\rangle}
  &\L{\langle p,p-2,p-3\rangle}
  &\L{\langle p-3,p\rangle}
  &\L{\langle (p-2)^{2},p-4\rangle}
  &\L{\langle p-3,p-4\rangle}\\\toprule
\langle (p-2)^{2}\rangle&.&.&1&1&.&.&1&.&.&.&.&.&.&.&.&.&.&1&.&.&.&.&.&.&.&.&.&.\\
\langle p-2,p-3\rangle&.&.&.&.&.&.&1&.&.&1&.&.&.&.&.&.&.&1&1&.&.&.&.&.&.&.&.&.\\
\langle p-3\rangle&.&.&.&1&.&.&.&.&.&.&.&.&.&.&.&.&.&1&1&1&.&.&.&.&.&.&.&.\\
\langle (p-2)^{2},p-1\rangle&1&1&1&1&1&.&1&1&.&.&.&.&.&.&.&.&.&1&.&.&1&.&.&.&.&.&.&.\\
\langle p-1,p-2,p-3\rangle&.&.&.&.&.&.&.&.&.&.&.&.&.&.&.&.&.&1&1&.&1&1&.&.&.&.&.&.\\
\langle p-3,p-1\rangle&.&1&.&1&.&.&1&1&.&1&.&.&.&.&.&.&.&1&1&1&1&1&1&.&.&.&.&.\\
\langle (p-2)^{2},p\rangle&1&.&.&.&.&1&.&1&1&.&.&.&.&.&.&.&.&.&.&.&1&.&.&1&.&.&.&.\\
\langle p,p-2,p-3\rangle&.&1&.&.&1&1&.&1&1&.&1&.&1&.&.&.&.&.&.&.&1&1&1&1&1&.&.&.\\
\langle p-3,p\rangle&.&1&.&.&.&.&.&.&.&.&.&.&.&.&.&.&.&.&.&.&1&1&.&1&1&1&.&.\\
\langle (p-2)^{2},p-4\rangle&.&.&.&.&.&1&.&.&1&.&1&1&.&1&1&1&1&.&.&.&.&.&.&1&1&.&1&.\\
\langle p-2,p-3,p-4\rangle&.&.&.&.&.&.&.&.&1&.&.&.&.&.&.&1&.&.&.&.&.&.&.&1&.&.&1&.\\
\langle p-3,p-4\rangle&.&.&.&.&.&1&.&.&.&.&.&.&1&1&.&.&1&.&.&.&.&.&.&1&1&1&1&1\\
\langle (p-2)^{3}\rangle&.&.&.&.&.&.&.&.&.&.&.&.&.&.&1&1&1&.&.&.&.&.&.&.&.&.&1&.\\
\langle (p-3)^{2},p-2\rangle&.&.&.&.&.&.&.&.&.&.&.&1&.&1&.&1&.&.&.&.&.&.&.&.&.&.&1&.\\
\langle (p-3)^{2}\rangle&.&.&.&.&.&.&.&.&.&.&.&.&.&1&.&.&1&.&.&.&.&.&.&.&.&.&1&1\\
\end{array}$
\end{sideways}

\subsection{Case 2: $(i,j)=(p-3, p-1)$}
\begin{sideways}
$\begin{array}{l|*{26}c}
  &\L{\langle p-3,p-2\rangle}
  &\L{\langle p,p-2\rangle}
  &\L{\langle p,p-2,p-3\rangle}
  &\L{\langle p-2,p-1\rangle}
  &\L{\langle p-3,p-1\rangle}
  &\L{\langle p-1,p-2\rangle}
  &\L{\langle p-1,p-2,p-3\rangle}
  &\L{\langle p,p-1,p-2\rangle}
  &\L{\langle p,p-1,p-3\rangle}
  &\L{\langle (p-2)^{2}\rangle}
  &\L{\langle (p-2)^{2},p-3\rangle}
  &\L{\langle (p-2)^{2},p\rangle}
  &\L{\langle (p-2)^{2},p-1\rangle}
  &\L{\langle p,p-2,p-4\rangle}
  &\L{\langle p-1,p-2,p-4\rangle}
  &\L{\langle p-4,p-3\rangle}
  &\L{\langle (p-3)^{2}\rangle}
  &\L{\langle p-4\rangle}
  &\L{\langle (p-3)^{2},p-2\rangle}
  &\L{\langle p-4,p-2\rangle}
  &\L{\langle (p-3)^{2},p\rangle}
  &\L{\langle p,p-3,p-4\rangle}
  &\L{\langle p-4,p\rangle}
  &\L{\langle (p-3)^{2},p-1\rangle}
  &\L{\langle p-1,p-3,p-4\rangle}
  &\L{\langle p-4,p-1\rangle}\\\toprule
\langle (p-3)^{2}\rangle&.&.&.&1&1&.&.&.&.&1&.&.&.&.&.&.&1&.&.&.&.&.&.&.&.&.\\
\langle p-3,p-4\rangle&.&.&.&.&.&.&.&.&.&1&.&.&.&.&.&.&1&.&.&.&.&.&.&.&.&.\\
\langle p-4\rangle&.&.&.&.&1&.&.&.&.&.&.&.&.&.&.&.&1&1&.&.&.&.&.&.&.&.\\
\langle (p-3)^{2},p-2\rangle&.&.&.&1&1&1&.&.&.&1&1&.&.&.&.&.&1&.&1&.&.&.&.&.&.&.\\
\langle p-2,p-3,p-4\rangle&.&.&.&.&.&.&.&.&.&.&.&.&.&.&.&.&1&.&1&.&.&.&.&.&.&.\\
\langle p-4,p-2\rangle&.&.&.&.&1&.&.&.&.&1&1&.&.&.&.&.&1&1&1&1&.&.&.&.&.&.\\
\langle (p-3)^{2},p\rangle&1&1&.&.&.&1&1&1&1&.&1&1&.&.&.&.&.&.&1&.&1&.&.&.&.&.\\
\langle p,p-3,p-4\rangle&1&.&1&.&.&.&1&.&.&.&1&1&.&1&.&1&.&.&1&1&1&1&.&.&.&.\\
\langle p-4,p\rangle&.&.&.&.&.&.&.&.&1&.&.&.&.&.&.&.&.&.&1&.&1&1&1&.&.&.\\
\langle (p-3)^{2},p-1\rangle&1&1&.&.&.&.&.&.&.&.&.&1&1&.&.&.&.&.&.&.&1&.&.&1&.&.\\
\langle p-1,p-3,p-4\rangle&.&.&1&.&.&.&1&1&1&.&.&1&1&1&1&.&.&.&.&.&1&1&.&1&1&.\\
\langle p-4,p-1\rangle&1&.&.&.&.&.&1&.&1&.&.&.&.&.&.&1&.&.&.&.&1&1&1&1&1&1\\
\langle (p-3)^{3}\rangle&.&.&.&.&.&.&.&.&.&.&.&.&1&.&1&.&.&.&.&.&.&.&.&1&1&.\\
\langle (p-4)^{2},p-3\rangle&.&.&.&.&.&.&.&.&.&.&.&.&1&.&.&.&.&.&.&.&.&.&.&1&.&.\\
\langle (p-4)^{2}\rangle&.&.&.&.&.&.&.&.&.&.&.&.&.&.&.&.&.&.&.&.&.&.&.&1&1&1\\
\end{array}$
\end{sideways}

\subsection{Case 2: $(i,j)=(p-2, p+1)$}
\begin{sideways}
$\begin{array}{l|*{31}c}
  &\L{\langle p-1,p\rangle}
  &\L{\langle p-1,p-2\rangle}
  &\L{\langle p-2,p\rangle}
  &\L{\langle p-2,p-1\rangle}
  &\L{\langle p,p-1,p-2\rangle}
  &\L{\langle p^{2}\rangle}
  &\L{\langle p^{2},p-1\rangle}
  &\L{\langle p^{2},p-2\rangle}
  &\L{\langle (p-1)^{2}\rangle}
  &\L{\langle (p-1)^{2},p\rangle}
  &\L{\langle (p-1)^{2},p-2\rangle}
  &\L{\langle p,p-3\rangle}
  &\L{\langle p-1,p-3\rangle}
  &\L{\langle p,p-1,p-3\rangle}
  &\L{\langle (p-1)^{2},p-3\rangle}
  &\L{\langle p-3,p-2\rangle}
  &\L{\langle (p-2)^{2},p-3\rangle}
  &\L{\langle p-4,p-1\rangle}
  &\L{\langle (p-1)^{2},p-4\rangle}
  &\L{\langle p-4,p-2\rangle}
  &\L{\langle (p-2)^{2}\rangle}
  &\L{\langle p-2,p-3\rangle}
  &\L{\langle p-3\rangle}
  &\L{\langle (p-2)^{2},p\rangle}
  &\L{\langle p,p-2,p-3\rangle}
  &\L{\langle p-3,p\rangle}
  &\L{\langle (p-2)^{2},p-1\rangle}
  &\L{\langle p-1,p-2,p-3\rangle}
  &\L{\langle p-3,p-1\rangle}
  &\L{\langle (p-2)^{2},p-4\rangle}
  &\L{\langle p-3,p-4\rangle}\\\toprule
\langle (p-2)^{2}\rangle&.&1&.&1&.&.&.&.&1&.&.&.&.&.&.&.&.&.&.&.&1&.&.&.&.&.&.&.&.&.&.\\
\langle p-2,p-3\rangle&.&.&.&.&.&.&.&.&1&.&.&.&1&.&.&.&.&.&.&.&1&1&.&.&.&.&.&.&.&.&.\\
\langle p-3\rangle&.&.&.&1&.&.&.&.&.&.&.&.&.&.&.&.&.&.&.&.&1&1&1&.&.&.&.&.&.&.&.\\
\langle (p-2)^{2},p\rangle&1&1&1&1&1&1&1&1&1&1&.&.&.&.&.&.&.&.&.&.&1&.&.&1&.&.&.&.&.&.&.\\
\langle p,p-2,p-3\rangle&.&.&.&.&.&1&.&.&1&1&.&1&1&1&.&.&.&.&.&.&1&1&.&1&1&.&.&.&.&.&.\\
\langle p-3,p\rangle&.&.&1&1&1&.&.&1&.&.&.&.&.&.&.&.&.&.&.&.&1&1&1&1&1&1&.&.&.&.&.\\
\langle (p-2)^{2},p-1\rangle&.&.&1&.&1&1&1&1&.&1&1&.&.&.&.&.&.&.&.&.&.&.&.&1&.&.&1&.&.&.&.\\
\langle p-1,p-2,p-3\rangle&.&.&.&.&1&1&.&.&.&.&.&1&.&.&.&.&.&.&.&.&.&.&.&1&1&.&1&1&.&.&.\\
\langle p-3,p-1\rangle&1&.&1&.&1&.&.&1&.&1&1&.&.&1&.&.&.&.&.&.&.&.&.&1&1&1&1&1&1&.&.\\
\langle (p-2)^{2},p-4\rangle&.&.&.&.&1&.&.&.&.&.&1&.&.&.&1&1&1&1&1&1&.&.&.&.&.&.&1&1&1&1&.\\
\langle p-2,p-3,p-4\rangle&.&.&.&.&1&.&.&.&.&.&1&.&.&.&.&.&.&.&1&.&.&.&.&.&.&.&1&.&.&1&.\\
\langle p-3,p-4\rangle&.&.&.&.&1&.&.&.&.&.&.&.&.&.&.&.&1&.&.&1&.&.&.&.&.&.&1&1&.&1&1\\
\langle (p-2)^{3}\rangle&.&.&.&.&1&.&.&.&.&.&.&.&.&.&.&.&.&1&1&1&.&.&.&.&.&.&.&.&.&1&.\\
\langle (p-3)^{2},p-2\rangle&.&.&.&.&.&.&.&.&.&.&.&.&.&.&1&.&1&.&1&.&.&.&.&.&.&.&.&.&.&1&.\\
\langle (p-3)^{2}\rangle&.&.&.&.&1&.&.&.&.&.&.&.&.&.&.&1&1&.&.&1&.&.&.&.&.&.&.&.&.&1&1\\
\end{array}$
\end{sideways}

\subsection{Case 2: $(i,j)=(p-3, p)$}
\begin{sideways}
$\begin{array}{l|*{23}c}
  &\L{\langle p-2,p-1\rangle}
  &\L{\langle p-1,p-2,p-3\rangle}
  &\L{\langle p-2,p\rangle}
  &\L{\langle p-3,p\rangle}
  &\L{\langle p,p-1,p-2\rangle}
  &\L{\langle p,p-2,p-3\rangle}
  &\L{\langle (p-1)^{2}\rangle}
  &\L{\langle (p-1)^{2},p-2\rangle}
  &\L{\langle (p-1)^{2},p-3\rangle}
  &\L{\langle (p-2)^{2}\rangle}
  &\L{\langle (p-2)^{2},p-1\rangle}
  &\L{\langle (p-2)^{2},p-3\rangle}
  &\L{\langle (p-2)^{2},p\rangle}
  &\L{\langle p,p-2,p-4\rangle}
  &\L{\langle p-4,p-3\rangle}
  &\L{\langle (p-3)^{2}\rangle}
  &\L{\langle p-4\rangle}
  &\L{\langle (p-3)^{2},p-1\rangle}
  &\L{\langle p-4,p-1\rangle}
  &\L{\langle (p-3)^{2},p-2\rangle}
  &\L{\langle p-4,p-2\rangle}
  &\L{\langle (p-3)^{2},p\rangle}
  &\L{\langle p,p-3,p-4\rangle}\\\toprule
\langle (p-3)^{2}\rangle&.&.&1&1&.&.&.&.&.&1&.&.&.&.&.&1&.&.&.&.&.&.&.\\
\langle p-3,p-4\rangle&.&.&.&.&.&.&.&.&.&1&.&.&.&.&.&1&.&.&.&.&.&.&.\\
\langle p-4\rangle&.&.&.&1&.&.&.&.&.&.&.&.&.&.&.&1&1&.&.&.&.&.&.\\
\langle (p-3)^{2},p-1\rangle&.&.&1&1&.&.&1&1&1&1&1&.&.&.&.&1&.&1&.&.&.&.&.\\
\langle p-1,p-3,p-4\rangle&.&.&.&.&.&.&1&.&.&1&1&.&.&.&.&1&.&1&.&.&.&.&.\\
\langle p-4,p-1\rangle&.&.&.&1&.&.&.&.&1&.&.&.&.&.&.&1&1&1&1&.&.&.&.\\
\langle (p-3)^{2},p-2\rangle&1&1&.&.&1&.&1&1&1&.&1&1&.&.&.&.&.&1&.&1&.&.&.\\
\langle p-2,p-3,p-4\rangle&.&.&.&.&.&.&1&.&.&.&.&.&.&.&.&.&.&1&.&1&.&.&.\\
\langle p-4,p-2\rangle&.&1&.&.&.&.&.&.&1&.&1&1&.&.&.&.&.&1&1&1&1&.&.\\
\langle (p-3)^{2},p\rangle&1&.&.&.&.&1&.&.&.&.&.&1&1&.&.&.&.&.&.&1&.&1&.\\
\langle p,p-3,p-4\rangle&.&1&.&.&1&1&.&.&.&.&.&1&1&1&1&.&.&.&.&1&1&1&1\\
\langle p-4,p\rangle&.&1&.&.&.&.&.&.&.&.&.&.&.&.&.&.&.&.&.&1&.&1&1\\
\langle (p-3)^{3}\rangle&.&.&.&.&.&1&.&.&.&.&.&.&1&1&.&.&.&.&.&.&.&1&1\\
\langle (p-4)^{2},p-3\rangle&.&.&.&.&.&.&.&.&.&.&.&.&1&.&.&.&.&.&.&.&.&1&.\\
\langle (p-4)^{2}\rangle&.&.&.&.&.&1&.&.&.&.&.&.&.&.&1&.&.&.&.&.&.&1&1\\
\end{array}$
\end{sideways}

\subsection{Case 2: $(i,j)=(p-3, p+1)$}
\begin{sideways}
$\begin{array}{l|*{23}c}
  &\L{\langle p-2,p-3\rangle}
  &\L{\langle p-3,p-2\rangle}
  &\L{\langle p,p-1,p-2\rangle}
  &\L{\langle p,p-1,p-3\rangle}
  &\L{\langle p,p-2,p-3\rangle}
  &\L{\langle p-1,p-2,p-3\rangle}
  &\L{\langle p^{2},p-2\rangle}
  &\L{\langle p^{2},p-3\rangle}
  &\L{\langle (p-1)^{2},p\rangle}
  &\L{\langle (p-1)^{2},p-2\rangle}
  &\L{\langle (p-1)^{2},p-3\rangle}
  &\L{\langle (p-2)^{2}\rangle}
  &\L{\langle (p-2)^{2},p\rangle}
  &\L{\langle (p-2)^{2},p-1\rangle}
  &\L{\langle (p-2)^{2},p-3\rangle}
  &\L{\langle (p-3)^{2}\rangle}
  &\L{\langle p-4\rangle}
  &\L{\langle (p-3)^{2},p\rangle}
  &\L{\langle p-4,p\rangle}
  &\L{\langle (p-3)^{2},p-1\rangle}
  &\L{\langle p-4,p-1\rangle}
  &\L{\langle (p-3)^{2},p-2\rangle}
  &\L{\langle p-4,p-2\rangle}\\\toprule
\langle (p-3)^{2}\rangle&1&1&.&.&.&.&.&.&.&.&.&1&.&.&.&1&.&.&.&.&.&.&.\\
\langle p-3,p-4\rangle&.&.&.&.&.&.&.&.&.&.&.&1&.&.&.&1&.&.&.&.&.&.&.\\
\langle p-4\rangle&.&1&.&.&.&.&.&.&.&.&.&.&.&.&.&1&1&.&.&.&.&.&.\\
\langle (p-3)^{2},p\rangle&1&1&.&.&1&.&1&1&.&.&.&1&1&.&.&1&.&1&.&.&.&.&.\\
\langle p,p-3,p-4\rangle&.&.&.&.&.&.&.&.&.&.&.&1&1&.&.&1&.&1&.&.&.&.&.\\
\langle p-4,p\rangle&.&1&.&.&1&.&.&1&.&.&.&.&.&.&.&1&1&1&1&.&.&.&.\\
\langle (p-3)^{2},p-1\rangle&.&.&1&1&1&1&1&1&1&1&1&.&1&1&.&.&.&1&.&1&.&.&.\\
\langle p-1,p-3,p-4\rangle&.&.&.&.&.&.&.&.&1&.&.&.&1&1&.&.&.&1&.&1&.&.&.\\
\langle p-4,p-1\rangle&.&.&.&1&1&1&.&1&.&.&1&.&.&.&.&.&.&1&1&1&1&.&.\\
\langle (p-3)^{2},p-2\rangle&.&.&.&1&.&1&.&.&1&1&1&.&.&1&1&.&.&.&.&1&.&1&.\\
\langle p-2,p-3,p-4\rangle&.&.&.&.&.&1&.&.&1&.&.&.&.&.&.&.&.&.&.&1&.&1&.\\
\langle p-4,p-2\rangle&.&.&1&1&.&1&.&.&.&.&1&.&.&1&1&.&.&.&.&1&1&1&1\\
\langle (p-3)^{3}\rangle&.&.&.&.&.&.&.&.&.&.&.&.&.&.&1&.&.&.&.&.&.&1&1\\
\langle (p-4)^{2},p-3\rangle&.&.&.&.&.&1&.&.&.&.&.&.&.&.&1&.&.&.&.&.&.&1&.\\
\langle (p-4)^{2}\rangle&.&.&.&.&.&.&.&.&.&.&.&.&.&.&.&.&.&.&.&.&.&1&.\\
\end{array}$
\end{sideways}

\subsection{Case 3: $(i,j,k)=(p-1, p+1, p+1)$}
$\begin{array}{l|*{12}c}
  &\L{\langle p\rangle}
  &\L{\langle p^{2}\rangle}
  &\L{\langle p,p-1\rangle}
  &\L{\langle p-1,p\rangle}
  &\L{\langle p^{2},p-1\rangle}
  &\L{\langle p^{3}\rangle}
  &\L{\langle p,p-2\rangle}
  &\L{\langle p^{2},p-2\rangle}
  &\L{\langle (p-1)^{2},p\rangle}
  &\L{\langle p,p-1,p-2\rangle}
  &\L{\langle (p-1)^{2}\rangle}
  &\L{\langle p-1,p-2\rangle}\\\toprule
\langle (p-1)^{2},p\rangle&1&1&1&.&.&1&.&.&1&.&.&.\\
\langle p,p-1,p-2\rangle&.&.&1&.&1&1&1&1&1&1&.&.\\
\langle (p-1)^{2}\rangle&.&1&1&.&1&1&.&.&1&.&1&.\\
\langle p-1,p-2\rangle&.&.&.&1&1&1&.&1&1&1&1&1\\
\end{array}$

\subsection{Case 3: $(i,j,k)=(p-2, p, p)$}
$\begin{array}{l|*{15}c}
  &\L{\langle p-1\rangle}
  &\L{\langle p-2,p-1\rangle}
  &\L{\langle (p-1)^{2},p-2\rangle}
  &\L{\langle p-1,p\rangle}
  &\L{\langle p,p-1\rangle}
  &\L{\langle (p-1)^{2},p\rangle}
  &\L{\langle p,p-1,p-2\rangle}
  &\L{\langle p-2,p\rangle}
  &\L{\langle (p-1)^{3}\rangle}
  &\L{\langle p-1,p-3\rangle}
  &\L{\langle (p-1)^{2},p-3\rangle}
  &\L{\langle (p-2)^{2},p-1\rangle}
  &\L{\langle p-1,p-2,p-3\rangle}
  &\L{\langle (p-2)^{2}\rangle}
  &\L{\langle p-2,p-3\rangle}\\\toprule
\langle (p-2)^{2},p-1\rangle&.&.&1&1&.&1&1&.&1&.&.&1&.&.&.\\
\langle p-1,p-2,p-3\rangle&.&.&.&1&.&.&.&.&1&1&1&1&1&.&.\\
\langle (p-2)^{2}\rangle&1&1&1&1&1&1&1&1&1&.&.&1&.&1&.\\
\langle p-2,p-3\rangle&.&.&.&.&.&.&.&.&1&.&1&1&1&1&1\\
\end{array}$

\subsection{Case 3: $(i,j,k)=(p-3, p-1, p-1)$}
$\begin{array}{l|*{8}c}
  &\L{\langle p-2,p-1\rangle}
  &\L{\langle p-1,p-2\rangle}
  &\L{\langle (p-2)^{2},p-1\rangle}
  &\L{\langle p-1,p-2,p-3\rangle}
  &\L{\langle p-3,p-1\rangle}
  &\L{\langle (p-2)^{3}\rangle}
  &\L{\langle (p-3)^{2},p-2\rangle}
  &\L{\langle (p-3)^{2}\rangle}\\\toprule
\langle (p-3)^{2},p-2\rangle&1&.&1&1&.&1&1&.\\
\langle p-2,p-3,p-4\rangle&1&.&.&.&.&1&1&.\\
\langle (p-3)^{2}\rangle&1&1&1&1&1&1&1&1\\
\langle p-3,p-4\rangle&.&.&.&.&.&1&1&1\\
\end{array}$

\subsection{Case 3: $(i,j,k)=(p-2, p, p+1)$}
$\begin{array}{l|*{16}c}
  &\L{\langle p-1,p-2\rangle}
  &\L{\langle p,p-1\rangle}
  &\L{\langle p-2,p-1\rangle}
  &\L{\langle (p-1)^{2},p\rangle}
  &\L{\langle p,p-1,p-2\rangle}
  &\L{\langle (p-1)^{2},p-2\rangle}
  &\L{\langle p-2,p\rangle}
  &\L{\langle p^{2},p-1\rangle}
  &\L{\langle (p-1)^{3}\rangle}
  &\L{\langle p-1,p-3\rangle}
  &\L{\langle p,p-1,p-3\rangle}
  &\L{\langle (p-1)^{2},p-3\rangle}
  &\L{\langle (p-2)^{2},p-1\rangle}
  &\L{\langle p-1,p-2,p-3\rangle}
  &\L{\langle (p-2)^{2}\rangle}
  &\L{\langle p-2,p-3\rangle}\\\toprule
\langle (p-2)^{2},p-1\rangle&.&.&.&1&.&.&.&.&1&.&.&.&1&.&.&.\\
\langle p-1,p-2,p-3\rangle&.&.&.&.&.&1&.&.&1&.&.&1&1&1&.&.\\
\langle (p-2)^{2}\rangle&1&1&1&1&1&1&1&1&1&.&.&.&1&.&1&.\\
\langle p-2,p-3\rangle&1&.&.&.&.&1&.&1&1&1&1&1&1&1&1&1\\
\end{array}$

\subsection{Case 3: $(i,j,k)=(p-3, p-1, p)$}
$\begin{array}{l|*{10}c}
  &\L{\langle (p-2)^{2},p-3\rangle}
  &\L{\langle p-2,p\rangle}
  &\L{\langle p,p-1,p-2\rangle}
  &\L{\langle (p-2)^{2},p\rangle}
  &\L{\langle p,p-2,p-3\rangle}
  &\L{\langle p-3,p\rangle}
  &\L{\langle (p-1)^{2},p-2\rangle}
  &\L{\langle (p-2)^{3}\rangle}
  &\L{\langle (p-3)^{2},p-2\rangle}
  &\L{\langle (p-3)^{2}\rangle}\\\toprule
\langle (p-3)^{2},p-2\rangle&1&.&.&1&1&.&.&1&1&.\\
\langle p-2,p-3,p-4\rangle&.&.&.&.&.&.&.&1&1&.\\
\langle (p-3)^{2}\rangle&1&1&1&1&1&1&1&1&1&1\\
\langle p-3,p-4\rangle&.&1&.&.&.&.&1&1&1&1\\
\end{array}$

\subsection{Case 3: $(i,j,k)=(p-3, p-1, p+1)$}
$\begin{array}{l|*{8}c}
  &\L{\langle (p-2)^{2},p-1\rangle}
  &\L{\langle p-1,p-2,p-3\rangle}
  &\L{\langle (p-2)^{2},p-3\rangle}
  &\L{\langle p-3,p-1\rangle}
  &\L{\langle (p-1)^{2},p-2\rangle}
  &\L{\langle (p-2)^{3}\rangle}
  &\L{\langle (p-3)^{2},p-2\rangle}
  &\L{\langle (p-3)^{2}\rangle}\\\toprule
\langle (p-3)^{2},p-2\rangle&1&.&.&.&.&1&1&.\\
\langle p-2,p-3,p-4\rangle&.&.&1&.&.&1&1&.\\
\langle (p-3)^{2}\rangle&1&1&1&1&1&1&1&1\\
\langle p-3,p-4\rangle&.&.&1&.&1&1&1&1\\
\end{array}$

\subsection{Case 3: $(i,j,k)=(p-2, p+1, p+1)$}
$\begin{array}{l|*{18}c}
  &\L{\langle p-1,p\rangle}
  &\L{\langle p^{2},p-1\rangle}
  &\L{\langle (p-1)^{2},p\rangle}
  &\L{\langle p,p-1,p-2\rangle}
  &\L{\langle p^{2},p-2\rangle}
  &\L{\langle p-2,p-1\rangle}
  &\L{\langle (p-1)^{2},p-2\rangle}
  &\L{\langle p^{3}\rangle}
  &\L{\langle (p-1)^{3}\rangle}
  &\L{\langle (p-2)^{2},p\rangle}
  &\L{\langle p,p-1,p-3\rangle}
  &\L{\langle p^{2},p-3\rangle}
  &\L{\langle p,p-2,p-3\rangle}
  &\L{\langle (p-1)^{2},p-3\rangle}
  &\L{\langle (p-2)^{2},p-1\rangle}
  &\L{\langle p-1,p-2,p-3\rangle}
  &\L{\langle (p-2)^{2}\rangle}
  &\L{\langle p-2,p-3\rangle}\\\toprule
\langle (p-2)^{2},p-1\rangle&1&1&1&1&.&.&.&1&1&1&.&.&.&.&1&.&.&.\\
\langle p-1,p-2,p-3\rangle&.&.&.&1&1&.&1&1&1&1&1&1&1&1&1&1&.&.\\
\langle (p-2)^{2}\rangle&.&.&1&1&.&.&1&.&1&.&.&.&.&.&1&.&1&.\\
\langle p-2,p-3\rangle&.&.&.&.&.&1&1&.&1&.&.&.&.&1&1&1&1&1\\
\end{array}$

\subsection{Case 3: $(i,j,k)=(p-3, p, p)$}
$\begin{array}{l|*{16}c}
  &\L{\langle p-2,p-1\rangle}
  &\L{\langle (p-1)^{2},p-3\rangle}
  &\L{\langle p-3,p-2\rangle}
  &\L{\langle (p-2)^{2},p-3\rangle}
  &\L{\langle p,p-1,p-2\rangle}
  &\L{\langle (p-1)^{2},p\rangle}
  &\L{\langle p,p-1,p-3\rangle}
  &\L{\langle p,p-2\rangle}
  &\L{\langle (p-2)^{2},p\rangle}
  &\L{\langle p,p-2,p-3\rangle}
  &\L{\langle p-3,p\rangle}
  &\L{\langle (p-1)^{3}\rangle}
  &\L{\langle (p-2)^{3}\rangle}
  &\L{\langle (p-3)^{2},p-1\rangle}
  &\L{\langle (p-3)^{2},p-2\rangle}
  &\L{\langle (p-3)^{2}\rangle}\\\toprule
\langle (p-3)^{2},p-2\rangle&.&1&.&1&1&1&1&.&1&1&.&1&1&1&1&.\\
\langle p-2,p-3,p-4\rangle&.&.&.&.&1&.&.&.&.&.&.&1&1&1&1&.\\
\langle (p-3)^{2}\rangle&1&.&1&1&1&.&.&1&1&1&1&.&1&.&1&1\\
\langle p-3,p-4\rangle&.&.&.&.&.&.&.&.&.&.&.&.&1&.&1&1\\
\end{array}$

\subsection{Case 3: $(i,j,k)=(p-3, p, p+1)$}
$\begin{array}{l|*{15}c}
  &\L{\langle p,p-1,p-2\rangle}
  &\L{\langle p-1,p-2,p-3\rangle}
  &\L{\langle (p-1)^{2},p\rangle}
  &\L{\langle (p-1)^{2},p-3\rangle}
  &\L{\langle p-3,p-2\rangle}
  &\L{\langle (p-2)^{2},p\rangle}
  &\L{\langle p,p-2,p-3\rangle}
  &\L{\langle (p-2)^{2},p-3\rangle}
  &\L{\langle p-3,p\rangle}
  &\L{\langle p^{2},p-2\rangle}
  &\L{\langle (p-1)^{3}\rangle}
  &\L{\langle (p-2)^{3}\rangle}
  &\L{\langle (p-3)^{2},p-1\rangle}
  &\L{\langle (p-3)^{2},p-2\rangle}
  &\L{\langle (p-3)^{2}\rangle}\\\toprule
\langle (p-3)^{2},p-2\rangle&.&.&1&.&.&1&.&.&.&.&1&1&1&1&.\\
\langle p-2,p-3,p-4\rangle&.&.&.&1&.&.&.&1&.&.&1&1&1&1&.\\
\langle (p-3)^{2}\rangle&1&1&.&.&1&1&1&1&1&1&.&1&.&1&1\\
\langle p-3,p-4\rangle&.&1&.&.&.&.&.&1&.&1&.&1&.&1&1\\
\end{array}$

\subsection{Case 3: $(i,j,k)=(p-3, p+1, p+1)$}
$\begin{array}{l|*{12}c}
  &\L{\langle p,p-1,p-2\rangle}
  &\L{\langle (p-1)^{2},p-2\rangle}
  &\L{\langle (p-2)^{2},p-1\rangle}
  &\L{\langle p-1,p-2,p-3\rangle}
  &\L{\langle (p-1)^{2},p-3\rangle}
  &\L{\langle p-3,p-2\rangle}
  &\L{\langle (p-2)^{2},p-3\rangle}
  &\L{\langle (p-1)^{3}\rangle}
  &\L{\langle (p-2)^{3}\rangle}
  &\L{\langle (p-3)^{2},p-1\rangle}
  &\L{\langle (p-3)^{2},p-2\rangle}
  &\L{\langle (p-3)^{2}\rangle}\\\toprule
\langle (p-3)^{2},p-2\rangle&1&1&1&1&.&.&.&1&1&1&1&.\\
\langle p-2,p-3,p-4\rangle&.&.&.&1&1&.&1&1&1&1&1&.\\
\langle (p-3)^{2}\rangle&.&.&1&1&.&.&1&.&1&.&1&1\\
\langle p-3,p-4\rangle&.&.&.&.&.&1&1&.&1&.&1&1\\
\end{array}$

\subsection{Case 4: $(i,j,k)=(p, p+1, p+1)$}
$\begin{array}{l|*{6}c}
  &\L{\langle p^{2}\rangle}
  &\L{\langle p,p-1\rangle}
  &\L{\langle p^{3}\rangle}
  &\L{\langle p^{2},p-1\rangle}
  &\L{\langle p-1,p\rangle}
  &\L{\langle p-1\rangle}\\\toprule
\langle p^{3}\rangle&.&.&1&.&.&.\\
\langle p^{2},p-1\rangle&1&.&1&1&.&.\\
\langle p-1,p\rangle&1&1&1&1&1&.\\
\langle p-1\rangle&.&.&1&1&1&1\\
\end{array}$

\subsection{Case 4: $(i,j,k)=(p-1, p, p)$}
$\begin{array}{l|*{11}c}
  &\L{\langle p-1\rangle}
  &\L{\langle (p-1)^{2}\rangle}
  &\L{\langle p-1,p\rangle}
  &\L{\langle p\rangle}
  &\L{\langle p,p-1\rangle}
  &\L{\langle (p-1)^{2},p\rangle}
  &\L{\langle p-1,p-2\rangle}
  &\L{\langle (p-1)^{3}\rangle}
  &\L{\langle (p-1)^{2},p-2\rangle}
  &\L{\langle p-2,p-1\rangle}
  &\L{\langle p-2\rangle}\\\toprule
\langle (p-1)^{3}\rangle&.&1&1&1&1&1&.&1&.&.&.\\
\langle (p-1)^{2},p-2\rangle&.&.&.&.&.&.&.&1&1&.&.\\
\langle p-2,p-1\rangle&1&1&1&.&.&1&1&1&1&1&.\\
\langle p-2\rangle&.&1&.&.&1&1&.&1&1&1&1\\
\end{array}$

\subsection{Case 4: $(i,j,k)=(p-2, p-1, p-1)$}
$\begin{array}{l|*{10}c}
  &\L{\langle p-2\rangle}
  &\L{\langle p-2,p-1\rangle}
  &\L{\langle p-1\rangle}
  &\L{\langle p-1,p-2\rangle}
  &\L{\langle (p-2)^{2},p-1\rangle}
  &\L{\langle p-2,p-3\rangle}
  &\L{\langle (p-2)^{3}\rangle}
  &\L{\langle (p-2)^{2},p-3\rangle}
  &\L{\langle p-3,p-2\rangle}
  &\L{\langle p-3\rangle}\\\toprule
\langle (p-2)^{3}\rangle&1&1&1&1&1&.&1&.&.&.\\
\langle (p-2)^{2},p-3\rangle&.&.&.&.&.&.&1&1&.&.\\
\langle p-3,p-2\rangle&.&1&.&.&1&1&1&1&1&.\\
\langle p-3\rangle&1&.&.&1&1&.&1&1&1&1\\
\end{array}$

\subsection{Case 4: $(i,j,k)=(p-3, p-2, p-2)$}
$\begin{array}{l|*{7}c}
  &\L{\langle p-3,p-2\rangle}
  &\L{\langle p-2\rangle}
  &\L{\langle p-2,p-3\rangle}
  &\L{\langle (p-3)^{2},p-2\rangle}
  &\L{\langle (p-3)^{3}\rangle}
  &\L{\langle p-4,p-3\rangle}
  &\L{\langle p-4\rangle}\\\toprule
\langle (p-3)^{3}\rangle&1&1&1&1&1&.&.\\
\langle (p-3)^{2},p-4\rangle&.&.&.&.&1&.&.\\
\langle p-4,p-3\rangle&1&.&.&1&1&1&.\\
\langle p-4\rangle&.&.&1&1&1&1&1\\
\end{array}$

\subsection{Case 4: $(i,j,k)=(p-1, p, p+1)$}
$\begin{array}{l|*{14}c}
  &\L{\langle p-1\rangle}
  &\L{\langle p\rangle}
  &\L{\langle (p-1)^{2}\rangle}
  &\L{\langle p,p-1\rangle}
  &\L{\langle (p-1)^{2},p\rangle}
  &\L{\langle p^{2}\rangle}
  &\L{\langle p^{2},p-1\rangle}
  &\L{\langle p-1,p-2\rangle}
  &\L{\langle p,p-2\rangle}
  &\L{\langle p,p-1,p-2\rangle}
  &\L{\langle (p-1)^{3}\rangle}
  &\L{\langle (p-1)^{2},p-2\rangle}
  &\L{\langle p-2,p-1\rangle}
  &\L{\langle p-2\rangle}\\\toprule
\langle (p-1)^{3}\rangle&1&1&1&.&.&1&1&.&.&.&1&.&.&.\\
\langle (p-1)^{2},p-2\rangle&.&.&1&1&1&1&1&1&1&1&1&1&.&.\\
\langle p-2,p-1\rangle&.&.&.&.&.&.&.&.&.&.&1&1&1&.\\
\langle p-2\rangle&1&.&1&.&1&.&1&1&.&1&1&1&1&1\\
\end{array}$

\subsection{Case 4: $(i,j,k)=(p-2, p-1, p)$}
$\begin{array}{l|*{15}c}
  &\L{\langle p-1,p-2\rangle}
  &\L{\langle (p-2)^{2},p-1\rangle}
  &\L{\langle p-2,p\rangle}
  &\L{\langle p-1,p\rangle}
  &\L{\langle p,p-1,p-2\rangle}
  &\L{\langle (p-2)^{2},p\rangle}
  &\L{\langle (p-1)^{2}\rangle}
  &\L{\langle (p-1)^{2},p-2\rangle}
  &\L{\langle p-2,p-3\rangle}
  &\L{\langle p-1,p-3\rangle}
  &\L{\langle p-1,p-2,p-3\rangle}
  &\L{\langle (p-2)^{3}\rangle}
  &\L{\langle (p-2)^{2},p-3\rangle}
  &\L{\langle p-3,p-2\rangle}
  &\L{\langle p-3\rangle}\\\toprule
\langle (p-2)^{3}\rangle&1&1&1&1&1&1&1&1&.&.&.&1&.&.&.\\
\langle (p-2)^{2},p-3\rangle&.&.&1&.&.&.&1&1&1&1&1&1&1&.&.\\
\langle p-3,p-2\rangle&.&1&.&.&.&1&.&.&.&.&.&1&1&1&.\\
\langle p-3\rangle&1&1&1&.&1&1&.&1&1&.&1&1&1&1&1\\
\end{array}$

\subsection{Case 4: $(i,j,k)=(p-3, p-2, p-1)$}
$\begin{array}{l|*{9}c}
  &\L{\langle p-3,p-1\rangle}
  &\L{\langle p-2,p-1\rangle}
  &\L{\langle p-1,p-2,p-3\rangle}
  &\L{\langle (p-3)^{2},p-1\rangle}
  &\L{\langle (p-2)^{2}\rangle}
  &\L{\langle (p-2)^{2},p-3\rangle}
  &\L{\langle (p-3)^{3}\rangle}
  &\L{\langle p-4,p-3\rangle}
  &\L{\langle p-4\rangle}\\\toprule
\langle (p-3)^{3}\rangle&1&1&1&1&1&1&1&.&.\\
\langle (p-3)^{2},p-4\rangle&1&.&.&.&1&1&1&.&.\\
\langle p-4,p-3\rangle&.&.&.&1&.&.&1&1&.\\
\langle p-4\rangle&1&.&1&1&.&1&1&1&1\\
\end{array}$

\subsection{Case 4: $(i,j,k)=(p-2, p-1, p+1)$}
$\begin{array}{l|*{12}c}
  &\L{\langle p-1,p-2\rangle}
  &\L{\langle p-1,p\rangle}
  &\L{\langle p,p-1,p-2\rangle}
  &\L{\langle (p-2)^{2},p-1\rangle}
  &\L{\langle (p-1)^{2}\rangle}
  &\L{\langle (p-1)^{2},p-2\rangle}
  &\L{\langle p-1,p-3\rangle}
  &\L{\langle p-1,p-2,p-3\rangle}
  &\L{\langle (p-2)^{3}\rangle}
  &\L{\langle (p-2)^{2},p-3\rangle}
  &\L{\langle p-3,p-2\rangle}
  &\L{\langle p-3\rangle}\\\toprule
\langle (p-2)^{3}\rangle&1&1&1&.&1&1&.&.&1&.&.&.\\
\langle (p-2)^{2},p-3\rangle&.&.&.&1&1&1&1&1&1&1&.&.\\
\langle p-3,p-2\rangle&.&.&.&.&.&.&.&.&1&1&1&.\\
\langle p-3\rangle&1&.&1&1&.&1&.&1&1&1&1&1\\
\end{array}$

\subsection{Case 4: $(i,j,k)=(p-3, p-2, p)$}
$\begin{array}{l|*{9}c}
  &\L{\langle (p-3)^{2},p-2\rangle}
  &\L{\langle p-2,p\rangle}
  &\L{\langle p,p-2,p-3\rangle}
  &\L{\langle (p-3)^{2},p\rangle}
  &\L{\langle (p-2)^{2}\rangle}
  &\L{\langle (p-2)^{2},p-3\rangle}
  &\L{\langle (p-3)^{3}\rangle}
  &\L{\langle p-4,p-3\rangle}
  &\L{\langle p-4\rangle}\\\toprule
\langle (p-3)^{3}\rangle&1&1&1&1&1&1&1&.&.\\
\langle (p-3)^{2},p-4\rangle&.&.&.&.&1&1&1&.&.\\
\langle p-4,p-3\rangle&1&.&.&1&.&.&1&1&.\\
\langle p-4\rangle&1&.&1&1&.&1&1&1&1\\
\end{array}$

\subsection{Case 4: $(i,j,k)=(p-3, p-2, p+1)$}
$\begin{array}{l|*{7}c}
  &\L{\langle p-2,p-1\rangle}
  &\L{\langle p-1,p-2,p-3\rangle}
  &\L{\langle (p-3)^{2},p-2\rangle}
  &\L{\langle (p-2)^{2}\rangle}
  &\L{\langle (p-2)^{2},p-3\rangle}
  &\L{\langle (p-3)^{3}\rangle}
  &\L{\langle p-4\rangle}\\\toprule
\langle (p-3)^{3}\rangle&1&1&.&1&1&1&.\\
\langle (p-3)^{2},p-4\rangle&.&.&1&1&1&1&.\\
\langle p-4,p-3\rangle&.&.&.&.&.&1&.\\
\langle p-4\rangle&.&1&1&.&1&1&1\\
\end{array}$

\subsection{Case 4: $(i,j,k)=(p-1, p+1, p+1)$}
$\begin{array}{l|*{13}c}
  &\L{\langle p\rangle}
  &\L{\langle p^{2}\rangle}
  &\L{\langle p^{2},p-1\rangle}
  &\L{\langle (p-1)^{2},p\rangle}
  &\L{\langle p^{3}\rangle}
  &\L{\langle p,p-2\rangle}
  &\L{\langle p,p-1,p-2\rangle}
  &\L{\langle p^{2},p-2\rangle}
  &\L{\langle p-2,p\rangle}
  &\L{\langle (p-1)^{3}\rangle}
  &\L{\langle (p-1)^{2},p-2\rangle}
  &\L{\langle p-2,p-1\rangle}
  &\L{\langle p-2\rangle}\\\toprule
\langle (p-1)^{3}\rangle&1&1&.&.&1&.&.&.&.&1&.&.&.\\
\langle (p-1)^{2},p-2\rangle&.&1&1&1&1&1&.&1&.&1&1&.&.\\
\langle p-2,p-1\rangle&.&.&.&1&1&.&1&1&1&1&1&1&.\\
\langle p-2\rangle&.&.&.&.&.&.&.&.&.&1&1&1&1\\
\end{array}$

\subsection{Case 4: $(i,j,k)=(p-2, p, p)$}
$\begin{array}{l|*{18}c}
  &\L{\langle p-2,p-1\rangle}
  &\L{\langle (p-1)^{2},p-2\rangle}
  &\L{\langle (p-2)^{2},p-1\rangle}
  &\L{\langle p-1,p\rangle}
  &\L{\langle p,p-1,p-2\rangle}
  &\L{\langle p,p-1\rangle}
  &\L{\langle (p-1)^{2},p\rangle}
  &\L{\langle p,p-2\rangle}
  &\L{\langle (p-2)^{2},p\rangle}
  &\L{\langle (p-1)^{3}\rangle}
  &\L{\langle p-1,p-3\rangle}
  &\L{\langle p-1,p-2,p-3\rangle}
  &\L{\langle (p-1)^{2},p-3\rangle}
  &\L{\langle p-3,p-1\rangle}
  &\L{\langle (p-2)^{3}\rangle}
  &\L{\langle (p-2)^{2},p-3\rangle}
  &\L{\langle p-3,p-2\rangle}
  &\L{\langle p-3\rangle}\\\toprule
\langle (p-2)^{3}\rangle&.&1&1&1&1&1&1&1&1&1&.&.&.&.&1&.&.&.\\
\langle (p-2)^{2},p-3\rangle&.&.&.&1&.&.&.&.&.&1&1&.&1&.&1&1&.&.\\
\langle p-3,p-2\rangle&1&1&1&.&1&.&1&.&1&1&.&1&1&1&1&1&1&.\\
\langle p-3\rangle&.&.&1&.&.&.&.&1&1&.&.&.&.&.&1&1&1&1\\
\end{array}$

\subsection{Case 4: $(i,j,k)=(p-3, p-1, p-1)$}
$\begin{array}{l|*{12}c}
  &\L{\langle p-3,p-2\rangle}
  &\L{\langle p-2,p-1\rangle}
  &\L{\langle p-1,p-2,p-3\rangle}
  &\L{\langle p-1,p-2\rangle}
  &\L{\langle (p-2)^{2},p-1\rangle}
  &\L{\langle p-1,p-3\rangle}
  &\L{\langle (p-3)^{2},p-1\rangle}
  &\L{\langle (p-2)^{3}\rangle}
  &\L{\langle p-4,p-2\rangle}
  &\L{\langle (p-3)^{3}\rangle}
  &\L{\langle p-4,p-3\rangle}
  &\L{\langle p-4\rangle}\\\toprule
\langle (p-3)^{3}\rangle&1&1&1&1&1&1&1&1&.&1&.&.\\
\langle (p-3)^{2},p-4\rangle&.&1&.&.&.&.&.&1&.&1&.&.\\
\langle p-4,p-3\rangle&.&.&1&.&1&.&1&1&1&1&1&.\\
\langle p-4\rangle&1&.&.&.&.&1&1&.&.&1&1&1\\
\end{array}$

\subsection{Case 4: $(i,j,k)=(p-2, p, p+1)$}
$\begin{array}{l|*{18}c}
  &\L{\langle p-2,p-1\rangle}
  &\L{\langle p,p-1\rangle}
  &\L{\langle p,p-1,p-2\rangle}
  &\L{\langle (p-2)^{2},p-1\rangle}
  &\L{\langle (p-1)^{2},p\rangle}
  &\L{\langle (p-2)^{2},p\rangle}
  &\L{\langle p^{2},p-1\rangle}
  &\L{\langle p^{2},p-2\rangle}
  &\L{\langle (p-1)^{3}\rangle}
  &\L{\langle p-1,p-2,p-3\rangle}
  &\L{\langle p,p-1,p-3\rangle}
  &\L{\langle (p-1)^{2},p-3\rangle}
  &\L{\langle p-3,p-1\rangle}
  &\L{\langle p,p-2,p-3\rangle}
  &\L{\langle (p-2)^{3}\rangle}
  &\L{\langle (p-2)^{2},p-3\rangle}
  &\L{\langle p-3,p-2\rangle}
  &\L{\langle p-3\rangle}\\\toprule
\langle (p-2)^{3}\rangle&1&1&.&1&.&.&1&1&1&.&.&.&.&.&1&.&.&.\\
\langle (p-2)^{2},p-3\rangle&.&.&1&1&1&1&1&1&1&1&1&1&.&1&1&1&.&.\\
\langle p-3,p-2\rangle&.&.&.&.&.&.&.&.&1&.&.&1&1&.&1&1&1&.\\
\langle p-3\rangle&1&.&.&1&.&1&.&1&.&1&.&.&.&1&1&1&1&1\\
\end{array}$

\subsection{Case 4: $(i,j,k)=(p-3, p-1, p)$}
$\begin{array}{l|*{15}c}
  &\L{\langle p-1,p-2,p-3\rangle}
  &\L{\langle (p-2)^{2},p-1\rangle}
  &\L{\langle (p-3)^{2},p-1\rangle}
  &\L{\langle p,p-2,p-3\rangle}
  &\L{\langle p,p-1,p-2\rangle}
  &\L{\langle (p-2)^{2},p\rangle}
  &\L{\langle p,p-1,p-3\rangle}
  &\L{\langle (p-3)^{2},p\rangle}
  &\L{\langle (p-1)^{2},p-2\rangle}
  &\L{\langle (p-1)^{2},p-3\rangle}
  &\L{\langle (p-2)^{3}\rangle}
  &\L{\langle p-4,p-2\rangle}
  &\L{\langle (p-3)^{3}\rangle}
  &\L{\langle p-4,p-3\rangle}
  &\L{\langle p-4\rangle}\\\toprule
\langle (p-3)^{3}\rangle&1&1&1&1&1&1&1&1&1&1&1&.&1&.&.\\
\langle (p-3)^{2},p-4\rangle&.&.&.&1&.&.&.&.&1&1&1&.&1&.&.\\
\langle p-4,p-3\rangle&.&1&1&.&.&1&.&1&.&.&1&1&1&1&.\\
\langle p-4\rangle&1&.&1&1&.&.&1&1&.&1&.&.&1&1&1\\
\end{array}$

\subsection{Case 4: $(i,j,k)=(p-3, p-1, p+1)$}
$\begin{array}{l|*{10}c}
  &\L{\langle p-1,p-2,p-3\rangle}
  &\L{\langle p,p-1,p-2\rangle}
  &\L{\langle (p-2)^{2},p-1\rangle}
  &\L{\langle p,p-1,p-3\rangle}
  &\L{\langle (p-3)^{2},p-1\rangle}
  &\L{\langle (p-1)^{2},p-2\rangle}
  &\L{\langle (p-1)^{2},p-3\rangle}
  &\L{\langle (p-2)^{3}\rangle}
  &\L{\langle (p-3)^{3}\rangle}
  &\L{\langle p-4\rangle}\\\toprule
\langle (p-3)^{3}\rangle&1&1&.&1&.&1&1&1&1&.\\
\langle (p-3)^{2},p-4\rangle&.&.&1&.&1&1&1&1&1&.\\
\langle p-4,p-3\rangle&.&.&.&.&.&.&.&1&1&.\\
\langle p-4\rangle&1&.&.&1&1&.&1&.&1&1\\
\end{array}$

\subsection{Case 4: $(i,j,k)=(p-2, p+1, p+1)$}
$\begin{array}{l|*{13}c}
  &\L{\langle p-1,p\rangle}
  &\L{\langle (p-1)^{2},p\rangle}
  &\L{\langle (p-1)^{2},p-2\rangle}
  &\L{\langle (p-2)^{2},p-1\rangle}
  &\L{\langle (p-1)^{3}\rangle}
  &\L{\langle p,p-1,p-3\rangle}
  &\L{\langle p-1,p-2,p-3\rangle}
  &\L{\langle (p-1)^{2},p-3\rangle}
  &\L{\langle p-3,p-1\rangle}
  &\L{\langle (p-2)^{3}\rangle}
  &\L{\langle (p-2)^{2},p-3\rangle}
  &\L{\langle p-3,p-2\rangle}
  &\L{\langle p-3\rangle}\\\toprule
\langle (p-2)^{3}\rangle&1&1&.&.&1&.&.&.&.&1&.&.&.\\
\langle (p-2)^{2},p-3\rangle&.&1&1&1&1&1&.&1&.&1&1&.&.\\
\langle p-3,p-2\rangle&.&.&.&1&1&.&1&1&1&1&1&1&.\\
\langle p-3\rangle&.&.&.&.&.&.&.&.&.&1&1&1&1\\
\end{array}$

\subsection{Case 4: $(i,j,k)=(p-3, p, p)$}
$\begin{array}{l|*{14}c}
  &\L{\langle p-1,p-2,p-3\rangle}
  &\L{\langle (p-2)^{2},p-3\rangle}
  &\L{\langle (p-3)^{2},p-2\rangle}
  &\L{\langle p,p-1,p-2\rangle}
  &\L{\langle p,p-2,p-3\rangle}
  &\L{\langle p,p-2\rangle}
  &\L{\langle (p-2)^{2},p\rangle}
  &\L{\langle p,p-3\rangle}
  &\L{\langle (p-3)^{2},p\rangle}
  &\L{\langle (p-2)^{3}\rangle}
  &\L{\langle p-4,p-2\rangle}
  &\L{\langle (p-3)^{3}\rangle}
  &\L{\langle p-4,p-3\rangle}
  &\L{\langle p-4\rangle}\\\toprule
\langle (p-3)^{3}\rangle&.&1&1&1&1&1&1&1&1&1&.&1&.&.\\
\langle (p-3)^{2},p-4\rangle&.&.&.&1&.&.&.&.&.&1&.&1&.&.\\
\langle p-4,p-3\rangle&1&1&1&.&1&.&1&.&1&1&1&1&1&.\\
\langle p-4\rangle&.&.&1&.&.&.&.&1&1&.&.&1&1&1\\
\end{array}$

\subsection{Case 4: $(i,j,k)=(p-3, p, p+1)$}
$\begin{array}{l|*{11}c}
  &\L{\langle p-1,p-2,p-3\rangle}
  &\L{\langle p,p-1,p-2\rangle}
  &\L{\langle p,p-2,p-3\rangle}
  &\L{\langle (p-3)^{2},p-2\rangle}
  &\L{\langle (p-2)^{2},p\rangle}
  &\L{\langle (p-3)^{2},p\rangle}
  &\L{\langle p^{2},p-2\rangle}
  &\L{\langle p^{2},p-3\rangle}
  &\L{\langle (p-2)^{3}\rangle}
  &\L{\langle (p-3)^{3}\rangle}
  &\L{\langle p-4\rangle}\\\toprule
\langle (p-3)^{3}\rangle&1&1&.&1&.&.&1&1&1&1&.\\
\langle (p-3)^{2},p-4\rangle&.&.&1&1&1&1&1&1&1&1&.\\
\langle p-4,p-3\rangle&.&.&.&.&.&.&.&.&1&1&.\\
\langle p-4\rangle&1&.&.&1&.&1&.&1&.&1&1\\
\end{array}$

\subsection{Case 4: $(i,j,k)=(p-3, p+1, p+1)$}
$\begin{array}{l|*{6}c}
  &\L{\langle p,p-1,p-2\rangle}
  &\L{\langle (p-2)^{2},p-1\rangle}
  &\L{\langle (p-2)^{2},p-3\rangle}
  &\L{\langle (p-3)^{2},p-2\rangle}
  &\L{\langle (p-2)^{3}\rangle}
  &\L{\langle (p-3)^{3}\rangle}\\\toprule
\langle (p-3)^{3}\rangle&1&1&.&.&1&1\\
\langle (p-3)^{2},p-4\rangle&.&1&1&1&1&1\\
\langle p-4,p-3\rangle&.&.&.&1&1&1\\
\langle p-4\rangle&.&.&.&.&.&1\\
\end{array}$

\bigskip
\[
\text{REFERENCES} 
\]

\bigskip

\begin{enumerate}
    \item[{[BO]}] \textsc{C. Bessenrodt and J.B. Olsson}, 
	On residue symbols and the Mullineux conjecture, \textit{J.
	Algebraic Comb.}, \textbf{7} (1998), 227-251.

    \item[{[D]}] \textsc{S. Donkin}, A note on decomposition numbers
	for general linear groups and symmetric groups, \textit{Math.
	Proc. Camb. Phil. Soc.}, \textbf{97 }(1985), 57-62.

    \item[{[BK]}] \textsc{J. Brundan and A. Kleshchev}, Representation
	theory of symmetric groups and their double covers,
	\textit{in} ``Proceedings of the Durham Conference on Groups,
	Representations, and Geometry'' (2000), \textit{to appear.}

    \item[{[FK]}] \textsc{B. Ford and A. Kleshchev}, A proof of the
	Mullineux conjecture, \textit{Math. Z.}, \textbf{226} (1997)),
	267-308.

    \item[{[J1]}] \textsc{G.D. James}, On the decomposition numbers of
	the symmetric groups III, \textit{J. Algebra}, \textbf{71}
	(1981), 115-122.

\item[{[J2]}] \leavevmode\vrule height 2pt depth -1.6pt width 23pt,
    Some combinatorial results involving Young diagrams, \textit{Math.
    Proc. Camb. Phil. Soc.}, \textbf{83} (1978), 1-10.

\item[{[J3]}] \leavevmode\vrule height 2pt depth -1.6pt width 23pt,
    The decomposition matrices of GL$_{n}(q)$ for $n\leq 10$,
    \textit{Proc. London Math. Soc.}, \textbf{60} (1990), 225-265.

\item[{[JK]}] \textsc{G.D. James and A. Kerber}, \textit{The
    representation theory of the symmetric group}, Encyclopedia of
    Mathematics and its Applications \textbf{16} (Cambridge University
    Press, 1981).

\item[{[JM1]}] \textsc{G.D. James and A. Mathas}, A $q$-analogue of
    the Jantzen-Schaper theorem, \textit{Proc. London Math. Soc.} (3),
    \textbf{74} (1997), 241-274.

\item[{[JM2]}] \leavevmode\vrule height 2pt depth -1.6pt width 23pt,
    Equating decomposition numbers for different primes, \textit{J.
    Algebra}, \textbf{258} (2002), 599-614.

\item[{[K]}] \textsc{A.S. Kleshchev}, Branching rules for symmetric
    groups and applications, \textit{in} ``Algebraic Groups and their
    Applications'' (R.W. Carter and J.  Saxl Eds.), 103-130,
    Kluwer, Dordrecht/Boston/ London, 1998.

\item[{[LLT]}] \textsc{A. Lascoux, B. Leclerc and J. Thibon}, Hecke algebras 
    at roots of unity and crystal bases of quantum affine algebras, 
    \textit{Comm. Math. Phys.}, \textbf{181} (1996),  205--263.

\item[{[LM]}] \textsc{F. L\"{u}beck and J. M\"{u}ller}, private communication.

\item[{[LN]}] \textsc{F. L\"{u}beck and M. Nuenh\"{o}ffer},
    Enumerating large orbits and direct condensation,
    \textit{Experiment. Math.}, \textbf{10} (2001), 197-205.

\item[{[M1]}] \textsc{A. Mathas}, Decomposition matrices of Hecke
    algebras of type $\bf A$, \textit{in} ``Gap: groups, algorithms and
    programming, 3.4.4'', M. Sch\"onert et al., RWTH Aachen, 1997.

\item[{[M2]}] \leavevmode\vrule height 2pt depth -1.6pt width 23pt,
    \textit{Hecke algebras and Schur algebras of the symmetric group},
    Univ. Lecture Notes, A.M.S., \textbf{15}, 1999.

\item[{[MR]}] \textsc{S. Martin and L. Russell}, Defect 3 blocks of
    symmetric group algebras, \textit{J. Algebra}, \textbf{213} (1999),
    304-339.

\item[{[Sch]}] \textsc{K.D. Schaper}, Charakterformeln f\"{u}r
    Weyl-Moduln und Specht-Moduln in Primcharacteristik,
    \textit{Diplomarbeit, Universit\"{a}t Bonn}, 1981.

\item[{[S1]}] \textsc{J.C. Scopes}, Cartan matrices and Morita
    equivalence for blocks of the symmetric groups, 
    \textit{J. Algebra}, \textbf{142} (1991), 441-455.

\item[{[S2]}] \leavevmode\vrule height 2pt depth -1.6pt width 23pt,
    Symmetric group blocks of defect two, \textit{Quart. J.  Math.
    Oxford} (2), \textbf{46} (1995), 201-234.

\item[{[R]}] \textsc{M. Richards}, Some decomposition numbers for
    Hecke algebras of general linear groups, \textit{Math. Proc. Camb.
    Phil. Soc.}, \textbf{119} (1996), 383-402.

\end{enumerate}
\end{document}